\documentclass[11pt, oneside]{article}

\usepackage{enumerate}
\usepackage{amscd}
\usepackage{amsmath}
\usepackage{epsfig}
\usepackage{amssymb}
\usepackage{amsthm}    
\usepackage{latexsym}
\usepackage{graphics}

\numberwithin{equation}{section}

\theoremstyle{plain}
\newtheorem{theorem}{Theorem}[section]
\newtheorem{proposition}[theorem]{Proposition} 
\newtheorem{lemma}[theorem]{Lemma} 
\newtheorem{corollary}[theorem]{Corollary} 
\newtheorem{conjecture}[theorem]{Conjecture} 
\theoremstyle{definition}
\newtheorem{definition}[theorem]{Definition} 
\newtheorem{question}[theorem]{Question} 
\theoremstyle{remark}
\newtheorem{remark}[theorem]{Remark} 
\newtheorem{example}[theorem]{Example} 

\newcommand{\sss}{\scriptscriptstyle}
\newcommand{\email}[1]{{\scriptsize{\it E-mail address}\/: {\rm #1}} }
\providecommand{\norm}[1]{\lVert#1\rVert}

\DeclareMathOperator{\diag}{diag}

\begin{document}

\begin{titlepage}
\title{Rigidity and Volume Preserving Deformation on Degenerate Simplices}

\author{Lizhao Zhang
\footnote{\email{lizhaozhang@alum.mit.edu}}
}
\date{}
\end{titlepage}

\maketitle


\begin{abstract}
Given a degenerate $(n+1)$-simplex in a $d$-dimensional space $M^d$
(Euclidean, spherical or hyperbolic space, and $d\geq n$),
for each $k$, $1\leq k\leq n$,
Radon's theorem induces a partition of the set of $k$-faces into two subsets.
We prove that if the vertices of the simplex vary smoothly in $M^d$ for $d=n$,
and the volumes of $k$-faces in one subset are constrained only to
decrease while in the other subset only to increase, then any sufficiently small motion
must preserve the volumes of all $k$-faces;
and this property still holds
in $M^d$ for $d\geq n+1$ if an invariant $c_{k-1}(\alpha^{k-1})$ of
the degenerate simplex has the desired sign.
This answers a question posed by the author,
and the proof relies on an invariant $c_k(\omega)$
we discovered for any $k$-stress
$\omega$ on a cell complex in $M^d$.
We introduce a characteristic polynomial of the degenerate simplex by defining
$f(x)=\sum_{i=0}^{n+1}(-1)^{i}c_i(\alpha^i)x^{n+1-i}$,
and prove that the roots of $f(x)$ are real for the Euclidean case.
Some evidence suggests
the same conjecture for the hyperbolic case.
\end{abstract}

\section{Introduction}
\label{section_introduction}

\subsection{Main results and motivations}
Let $M^d$ of dimension $d\geq n$
be the Euclidean, spherical or hyperbolic space of
constant curvature $\kappa$,
and $\mathbf{A}=\{A_1, \dots, A_{n+2}\}$ be a set of
vertices of a degenerate $(n+1)$-dimensional simplex in $M^d$,
where by degenerate we mean the vertices are confined in a lower dimensional $M^n$.
Assume further that all $n$-faces of $\mathbf{A}$ are non-degenerate.
By Radon's theorem the vertices of $\mathbf{A}$ can be partitioned into two subsets whose
convex hulls in $M^d$ intersect.
The only trivial exception is for the spherical case
when the vertices are not confined in any open half sphere,
then in this case one subset of vertices should be the empty set.
For each $k$, $1\leq k\leq n$,
counting each $k$-face's number of vertices mod $2$ in each subset
induces a partition of the set of $k$-faces into two subsets
$X_{1,k}$ and $X_{2,k}$.
The author asked the following question in \cite{Zh}:

\begin{question}
\label{question_volume_preserve}
If $\mathbf{A}$ varies smoothly in $M^d$, and the volumes of
 $k$-faces in one subset ($X_{1,k}$ \emph{or} $X_{2,k}$)
are constrained only to
decrease while in the other subset only to increase,
does the motion preserve the volumes of all $k$-faces of $\mathbf{A}$?
\end{question}

The purpose of this paper is twofold. First,
we prove a rigidity theorem which
gives an affirmative answer to
Question \ref{question_volume_preserve} for $d=n$,
and shows that it still holds
for $d\geq n+1$ if an invariant $c_{k-1}(\alpha^{k-1})$ we obtained from
$\mathbf{A}$ has the desired sign.
Second, under the motivation of Question \ref{question_volume_preserve},
we develop a theory to link $k$-stress 
(a notion introduced by Lee \cite{Lee}, see also \cite{Ryb,TWW})
with the volume deformation on cell complexes
(not necessarily simplicial) in $M^d$,
discover a geometric invariant $c_k(\omega)$
for any $k$-stress
$\omega$ on a cell complex in $M^d$,
and introduce a notion of characteristic polynomial
of a degenerate simplex,
which is also of interest by its own right.
These two topics are strongly related.
To some extent, Question \ref{question_volume_preserve}
serves the purpose of storytelling,
which leads to the development of the theory of the second topic above.

To state our results, we first introduce some basic notions.
Let the spherical space $\mathbb{S}^d$ be the standard unit sphere
centered at the origin in a Euclidean space $\mathbb{R}^{d+1}$,
and the hyperbolic space $\mathbb{H}^d$ be described by the \emph{hyperboloid model}:
Let $\mathbb{R}^{d,1}$ be a $(d+1)$-dimensional vector space endowed with a
metric
\[ x\cdot y
=-x_0y_0+x_1y_1+\cdots+x_dy_d,
\]
then $\mathbb{H}^d$ is defined by
\[ \{x \in \mathbb{R}^{d,1}: x\cdot x = -1, \quad x_0 > 0\},
\]
which is the upper sheet of a two-sheeted hyperboloid.
Under this embedding, we can use the vector space to
discuss the \emph{linear} relations
between points in $\mathbb{S}^d$ or $\mathbb{H}^d$.

Since every $n$-face of $\mathbf{A}$ is non-degenerate,
so up to a constant factor,
there is an unique affine dependence among the vertices of $\mathbf{A}$
for the Euclidean case,
or a linear dependence for the non-Euclidean case.
Namely, there is a sequence of non-zero coefficients
$\alpha_{1}$, \ldots, $\alpha_{n+2}\in \mathbb{R}$,
such that
\begin{equation}
\label{equation_alpha_r_s_h}
\begin{split}
&\sum\alpha_{i}A_i=0 \quad\text{and}\quad \sum\alpha_{i}=0 \quad\text{(for the Euclidean case)},  \\
&\sum\alpha_{i}A_i=0 \quad\text{(for the spherical or hyperbolic case)}.  \\
\end{split}
\end{equation}
We call $\alpha:=\{\alpha_1, \dots, \alpha_{n+2}\}$
a \emph{$1$-stress} on $\mathbf{A}$.
We reserve the notations $\mathbf{A}$ and $\alpha$,
or simply $(\mathbf{A},\alpha)$,
as well as $G_{n,k}$ and $G'_{n,k}$ defined next,
for the rest of this paper.

\begin{definition}
\label{definition_framework}
Let $(\mathbf{A}, \alpha)$ be as in (\ref{equation_alpha_r_s_h})
where $\alpha$ is a 1-stress on $\mathbf{A}$.
For each $k$, $1\leq k\leq n$, define $G_{n,k}$ to be a framework
equipped with the following volume constraints on $k$-faces of $\mathbf{A}$:
the volume of a $k$-face $F$ is constrained only to decrease (under \emph{tension})
if $\prod_{A_s\in F}\alpha_s <0$,
and only to increase (under \emph{compression})
if $\prod_{A_s\in F}\alpha_s >0$.
And define $G'_{n,k}$ by flipping the tension-compression
volume constraints in $G_{n,k}$.
\end{definition}


Let $\mathbf{A}(t)$ be a smooth motion of $\mathbf{A}$ in $M^d$,
and $\mathbf{A}(0)=\mathbf{A}$ be the initial position.
Then our rigidity theorem can be formulated as follows.

\begin{theorem}
\emph{(Main Theorem 1)}
\label{theorem_mainone_r_s_h_n}
If $\mathbf{A}(t)$ varies smoothly over $t$ in $M^n$, 
then for both $G_{n,k}$ and $G'_{n,k}$
that equipped with the volume constraints on $\mathbf{A}$,
the motion must preserve the volumes of all k-faces of $\mathbf{A}(t)$
for small $t\geq 0$.
\end{theorem}

The case $d\geq n+1$ is much harder and very different
from the case $d=n$. One of the most important results
of this paper is an invariant $c_k(\omega)$
(Theorem \ref{theorem_invariant_r_s_h})
we obtained from any $k$-stress $\omega$
on a cell complex in $M^d$.
Particularly for $\mathbf{A}$
we derive a sequence of invariants
$c_0(\alpha^0), \dots, c_{n+1}(\alpha^{n+1})$ (Definition \ref{defintion_invariant_c_k}),
which plays a key role in both the formulation and proof
of the following theorem.

\begin{theorem}
\emph{(Main Theorem 2)}
\label{theorem_mainone_r_s_h_n_plus_1}
For $d\geq n+1$, if $\mathbf{A}(t)$ varies smoothly over $t$ in $M^d$
and $c_{k-1}(\alpha^{k-1})> 0$ (resp. $c_{k-1}(\alpha^{k-1})< 0$),
then for $G_{n,k}$ (resp. $G'_{n,k}$)
that equipped with the volume constraints on $\mathbf{A}$,
the motion must preserve the volumes of all k-faces of $\mathbf{A}(t)$,
and the vertices are confined in a lower dimensional $M^n$
for small $t\geq 0$.
\end{theorem}

For $k=n$ with $n\geq 2$, the statement that 
``the vertices are confined in a lower dimensional $M^n$ for small $t\geq 0$''
is somewhat surprising, 
because the number of volume constraints $n+2$
is far less than the degree
of freedom of $\mathbf{A}$ in $M^d$ up to congruence,
which is $(n+2)(n+1)/2$, or subtract by 1 if
$\mathbf{A}$ is restricted in $M^n$.
Note that in both Theorem \ref{theorem_mainone_r_s_h_n}
and Theorem \ref{theorem_mainone_r_s_h_n_plus_1},
except for $k=1$
we do \emph{not} prove that the motion is \emph{rigid},
which is a stronger notion than the type of volume rigidity we proved.
This can be a potential improvement to our results,
and will be addressed 
in Section \ref{section_related_questions} along with
some related questions.

As remarked above,
a key tool we use to prove the rigidity theorem is $k$-stress,
a notion first introduced by Lee 
on simplicial complexes with vertices chosen in the Euclidean space \cite{Lee}.
The introduction of the notion
was partly inspired by Kalai's proof of the
lower bound theorem using classical stresses \cite{Ka},
and motivated to give a geometric understanding of 
Stanley's proof of the necessity of the $g$-theorem
for simplicial convex polytopes \cite{St}, which used algebraic geometry.
A notable property of $k$-stress is that
for a simplicial $(d-1)$-sphere $\Delta$ with vertices chosen generically in 
$\mathbb{S}^{d-1}$, 
according to Lee \cite{Lee},
the dimension of the space of $k$-stresses on $\Delta$ in $\mathbb{S}^{d-1}$
is $h_k$, where $(h_0,\dots,h_d)$ is the $h$-vector of $\Delta$.
What remains open, which if true proves the $g$-conjecture for
simplicial spheres,
is to show that for $\Delta$ with vertices chosen generically
in $\mathbb{R}^d$,
the dimension of the space of $k$-stresses on $\Delta$ in $\mathbb{R}^d$
is $g_k$ for $k\leq\lfloor d/2\rfloor$, where $g_0=h_0$,
$g_k=h_k-h_{k-1}$, $k=1,\dots,\lfloor d/2\rfloor$.
As $k$-stresses are the central theme of
this paper and
non-degenerate simplices do not admit
any $k$-stress,
that is why the $(n+1)$-simplices we looked
at in this paper are degenerate.

For $k=1$, where $G_{n,1}$ and $G'_{n,1}$ 
are \emph{tensegrity frameworks},
Bezdek and Connelly proved
the Euclidean case of Theorem \ref{theorem_mainone_r_s_h_n}
and \ref{theorem_mainone_r_s_h_n_plus_1} \cite{BC}.
They actually proved a stronger version:
$G_{n,1}$ is \emph{globally rigid} in $\mathbb{R}^d$
for any $d\ge n$.
Rigidity and flexibility of tensegrity frameworks,
which analyze geometric structures equipped
with distance constraints on edges,
have been extensively investigated in the past (see \cite{Co2}). 
However, the analogue for $k\geq 2$
with \emph{volume constraints} on $k$-faces,
especially for the non-Euclidean case,
has been much less studied in the literature.

For $k=n$, recall the notation $X_{1,n}$ and $X_{2,n}$ 
as in Question \ref{question_volume_preserve}.
Note that in $\mathbb{R}^n$ both
$G_{n,n}$ and $G'_{n,n}$ must preserve the volumes of all $n$-faces, 
as under a continuous motion
the sum of the volumes of $n$-faces in $X_{1,n}$
is equal to those in $X_{2,n}$.
Also note that $G_{n,n}$ and $G'_{n,n}$ are not rigid in $\mathbb{R}^n$, 
as under affine motions they can change the shapes smoothly 
while preserving the volumes of all $n$-faces.
However, it is far from trivial to tell
if $G_{n,n}$ or $G'_{n,n}$ will still preserve the volumes of
$n$-faces in $\mathbb{R}^{n+1}$,
as potentially the vertices can be \emph{lifted} in $\mathbb{R}^{n+1}$
to form a non-degenerate $(n+1)$-simplex,
therefore the sum of the volumes
in $X_{1,n}$ is no longer necessarily equal to those in $X_{2,n}$.
For $n=2$, from our results
we come up with a particularly interesting example of
``four points on a circle'' to address this phenomenon
(Example \ref{example_four_points_2}),
which we present here as well:

In $\mathbb{R}^3$, given four points that are initially in convex position
in a 2-dimensional plane.
If we allow the four points to move smoothly in $\mathbb{R}^3$
but constrain all the triangles formed by any three points to 
preserve the areas during the motion,
then in order for the four points to form a \emph{non}-degenerate 3-simplex
in $\mathbb{R}^3$, they have to be confined in a plane first
until they move on to a common circle.
And only from this common circle they can be lifted to form a 
non-degenerate 3-simplex.%
\footnote{
It will be interesting to see if this phenomenon can be demonstrated in the ``real'' world
by using some physical material,
e.g., just as the minimal surface can be visualized by using soap film.
}

A similar analogue for the non-Euclidean case is also given in
Example \ref{example_four_points_2}.
In fact, these examples were part of the
motivations for the author to pose Question \ref{question_volume_preserve}
and formulate the rigidity theorem in the first place.

\subsection{Strategy overview}
Our strategy to prove the rigidity theorem is as follows.
Using the Schl\"{a}fli differential formula, we develop techniques
for $k$-stresses on cell complexes in $M^d$
(Theorem \ref{theorem_stress_volume} and \ref{theorem_invariant_r_s_h}).
Applying them 
we obtain a differential \emph{equality}
(Proposition \ref{proposition_stress_volume_A})
for the $k$-faces of $\mathbf{A}$ if $d=n$,
and a differential \emph{inequality}
(Proposition \ref{proposition_inequality_r_s_h}) if $d\geq n+1$,
which directly lead to the proofs of
Theorem \ref{theorem_mainone_r_s_h_n} 
and \ref{theorem_mainone_r_s_h_n_plus_1}
respectively.
We also obtain a new version of Schl\"{a}fli differential formula on simplices based
on edge lengths (Proposition \ref{proposition_schlafli_edge}).
Some remarks on the history of the Schl\"{a}fli differential formula can be found
in Milnor's paper \cite{Mil}.

To analyze the interrelation between the rigidity properties
of different dimensions $k$,
we introduce a notion of \emph{characteristic polynomial} of the
degenerate $(n+1)$-simplex by defining
$f(x)=\sum_{i=0}^{n+1}(-1)^{i}c_i(\alpha^i)x^{n+1-i}$.
For the Euclidean case, we prove that the roots of $f(x)$
are real and give a way to count the number of positive roots
(Theorem \ref{theorem_realroots_r}).
Some evidence suggests the same conjecture for the hyperbolic
case (Conjecture \ref{conjecture_realroots_h}).
And in Section \ref{section_char_poly_general}, we naturally generalize
the notion of characteristic polynomial $f(x)$ to a set of points
(continuous distribution allowed)
in $M^n$ associated with a $1$-stress
on the points.

\subsection{Historical works}
Rigidity and deformation of geometric structures
have attracted the attention of mathematicians for a long time.
One of the first substantial mathematical results concerning
rigidity is Cauchy's rigidity theorem,
which proved that all convex polyhedra with solid faces 
and flexible dihedral angles are rigid.
It was widely believed and conjectured that the same held true
for non-convex polyhedron as well.
However, Connelly disproved the rigidity conjecture
by constructing a flexible polyhedron in $\mathbb{R}^3$ \cite{Co1},
and with D. Sullivan, they conjectured that the volume
bounded by a flexible polyhedron is constant during the flex.
Sabitov proved the conjecture of Connelly and Sullivan
for flexible polyhedron homeomorphic to a sphere \cite{Sa};
and Connelly, Sabitov, and Walz proved it
for general polyhedral surface
in ``The bellows conjecture'' \cite{CSW}.
The same conjecture in the
spherical space is not true though. Alexandrov constructed
a flexible polyhedron in an open half sphere in $\mathbb{S}^3$ which does
not conserve the volume \cite{Al}.

Motivated by the historical works on rigidity, we bring a different
view to the field. Instead of analyzing geometric structures
equipped with distance constraints between vertices
(as in the above works), 
we analyze volume constraints on $k$-faces of the underlying
geometric structure, as well as the interrelation
between the rigidity properties of different dimensions $k$.
We also generalize our main results to the Euclidean,
spherical and hyperbolic space together, so the validity of
our rigidity theorem
is independent of the constant curvature \emph{value} of the underlying
space.

\section{Volume preserving deformation}
\label{section_r_s_h}

To prove the main results,
our approach emphasizes on the non-Euclidean case,
and treats the Euclidean case as a limit of the spherical case.

\subsection{Basic terminology}

The following terminology and definitions are intended to clarify
the meaning of terms used in this paper.
Some terms are new.

By a $k$-dimensional \emph{convex polytope} in $M^k$
we mean a compact subset which can be
expressed as a finite intersection of closed half spaces.
A \emph{cell complex} in $M^d$ is a finite set 
of convex polytopes (called \emph{cells}) in $M^d$,
such that every face (empty set included)
of a cell is also a cell in the set,
and any two cells share a unique maximal common face, the \emph{intersection}.
However, in this paper we do not worry about the self-intersections
between the cells in $M^d$.
For the spherical case, we also require that
each cell of a cell complex lies strictly in an open half sphere,
so $\mathbb{S}^0$ is not a cell and a 0-cell always contains only one point.
Also a half circle is not a cell.

We want to point out that the \emph{convexity} of the polytopes above
plays almost no role in the context of this paper.
However, for simplicity, we content ourselves with only considering
convex polytopes in $M^d$.

For a cell complex in $M^d$,
we call it a \emph{k-tensegrity framework} if it is
equipped with volume constraints on $k$-faces (equalities and inequalities, 
as tension and compression),
e.g., $G_{n,k}$ and $G'_{n,k}$.
A $k$-tensegrity framework $p$ is \emph{rigid} in $M^{d}$,
if any continuous motion in $M^d$ that satisfies the volume constraints
is also a rigid motion;
and it is \emph{globally rigid} in $M^{d}$,
if for any other configuration $q$ in $M^d$ satisfying the volume constraints,
$q$ is congruent to $p$.

The notion of $k$-tensegrity framework,
a new term introduced in this paper and \cite{Zh},
is a natural higher dimensional generalization
of the notion of \emph{tensegrity frameworks} (see \cite{Co2}),
which is a finite graph
with vertices in $M^d$ and
equipped with length constraints on edges.

\subsection{Stresses on cell complex}
\label{section_stresses}

The notion of $k$-stresses on cell complexes
in $M^d$ plays an important role
in proving the main rigidity theorem.
While our rigidity theorem concerns the boundary complex
of a degenerate simplex in $M^d$,
our results about $k$-stresses are much more general,
which can be extended to
cell complexes (not necessarily simplicial) in $M^d$
without much extra effort.

If $K$ is a cell complex in $M^d$,
with a slight abuse of notation we simply denote by $K$ as well the set of all its cells,
and by $K^r$ the subset of its $r$-cells.

\begin{definition}
\label{defintion_k_stress}
Consider a cell complex $K$ (not necessarily
of dimension $d-1$ or $d$) in $M^d$.
A \emph{$k$-stress} $(2\leq k\leq d+1)$ on $K$ is
a real-valued function $\omega$ on the $(k-1)$-cells of $K$,
such that for each $(k-2)$-cell $F$ of $K$,
\[\sum_{G\in K^{k-1},F\subset G}\omega(G)u_{\sss F,G}=0,
\]
where the sum is taken over all $(k-1)$-cells $G$ of $K$
that contain $F$,
and $u_{\sss F,G}$ is the inward unit normal to $G$ at its facet $F$.
For $k=1$, a $1$-stress is an affine dependence among the vertices
for the Euclidean case, or a linear dependence for the non-Euclidean case.
\end{definition}

The notion of $k$-stress
was first introduced by Lee \cite{Lee}
on simplicial complexes with vertices chosen in the Euclidean space
with a slightly different setting.
Lee considered two types of $k$-stresses, affine and linear.
For a simplicial complex $K$ with vertices chosen in $\mathbb{R}^d$,
the space of \emph{affine} $k$-stresses is isomorphic to
the space of our notion of $k$-stresses.
Denote by $b_a$ (resp. $b_l$) the affine (resp. linear) $k$-stress
on $K$ in $\mathbb{R}^d$,
then $\omega(G)=(k-1)!V_{k-1}(G)b_a(G)$ for each $(k-1)$-face $G$ of $K$,
where $\omega$ is our notion of $k$-stress
and $V_{k-1}(G)$ denotes the $(k-1)$-dimensional volume of $G$.
If $K$ is a spherical simplicial complex with vertices chosen
in $\mathbb{S}^{d-1}$, as $\mathbb{S}^{d-1}$ is embedded in $\mathbb{R}^d$, 
we can also loosely treat $K$ as a \emph{Euclidean} simplicial complex
 in $\mathbb{R}^d$ in the sense of Lee.
Under this interpretation,
the space of our notion of $k$-stresses on $K$ in $\mathbb{S}^{d-1}$
is isomorphic to the space of \emph{linear} $k$-stresses on $K$ in $\mathbb{R}^d$.
For a $(k-1)$-face $G$ of $K$ in $\mathbb{S}^{d-1}$,
let $\norm{G}$ be $k!$ times the volume of the Euclidean $k$-simplex
formed by the vertices of $G$ and the origin $O$, 
then $\omega(G)=\norm{G}b_l(G)$. Note that unlike $\omega$,
the linear $k$-stress $b_l$ cannot be extended to non-simplicial
cell complexes in $\mathbb{S}^{d-1}$, as $\norm{G}$ cannot be properly defined for
spherical cells $G$ that are not simplicial.

Rybnikov in \cite{Ryb} extended
the notion of $k$-stress to cell complexes in Euclidean and spherical spaces,
and our terminology agrees with its terminology.
Similar notions were also considered in \cite{TWW}.
McMullen also considered \emph{weights} on simple polytopes \cite{Mc1},
a notion dual to stresses.
The relationship between $k$-stresses and volumes of simplicial or 
simple polytopes in the Euclidean case
was discussed in \cite{Lee} and \cite{Mc1}. However, it seems
that our work is the first to give a systematic discussion of the relationship
between $k$-stresses and the volumes of faces of cell complexes in
the \emph{non-Euclidean} case.

\subsection{A differential formula}

As a first step to proving Theorem \ref{theorem_mainone_r_s_h_n},
we obtain a differential formula for
$(k+1)$-stresses in Theorem \ref{theorem_stress_volume},
which also establishes a correspondence between the signs 
of volume constraints and the signs of $(k+1)$-stresses on $k$-faces
of a cell complex.
It generalizes the well established correspondence 
of classical stresses on 1-dimensional faces of a framework.

For a $k$-polytope $G$ in $M^d$,
denote by $V_k(G)$ the $k$-dimensional volume of $G$.
To compute the differential of the volumes
of $k$-dimensional polytopes in $\mathbb{S}^d$ or $\mathbb{H}^d$,
Schl\"{a}fli's differential formula plays a central role.
Some remarks on the history of the Schl\"{a}fli differential formula can be
found in Milnor's paper \cite{Mil}.
Consider a family of $k$-dimensional convex polytopes $P$
which vary smoothly
in a space of constant curvature $\kappa$.
For each $(k-2)$-dimensional face $F$
let $\theta_{F}$ be the dihedral angle at $F$.
Then the Schl\"{a}fli differential formula states that
\begin{equation}
\label{equation_schlafli}
\kappa\cdot dV_k(P)=\frac{1}{k-1}\sum_{F}V_{k-2}(F)\,d\theta_{F},
\end{equation}
where the sum is taken over all $(k-2)$-faces $F$ of $P$.
When $k-2=0$, $V_0(F)$ is the number of points in $F$.

For each $(k-2)$-face $F$ of $P$, it can be uniquely
described as an intersection $F=E\cap E'$
of two $(k-1)$-faces $E$ and $E'$ of $P$.
Let $u_{\sss E,P}$ be the inward unit normal to $P$ at its facet $E$,
$u_{\sss F,E}$ be the inward unit normal to $E$ at its facet $F$, and so on.
Note that $u_{\sss E,P}$, $u_{\sss E',P}$, $u_{\sss F,E}$ and $u_{\sss F,E'}$
are all in a single $2$-dimensional plane;
the angle between $u_{\sss E,P}$ and $u_{\sss E',P}$ is $\pi-\theta_F$;
the angle between  $u_{\sss F,E}$ and $u_{\sss F,E'}$ is $\theta_F$;
the angle between $u_{\sss E,P}$ and $u_{\sss F,E}$ is $\pi/2$;
and the angle between $u_{\sss E',P}$ and $u_{\sss F,E'}$ is $\pi/2$ as well.
It is easy to check that
\begin{equation}
\label{equation_differential_theta_F}
d\theta_F=-u_{\sss E,P}\cdot du_{\sss F,E}-u_{\sss E',P}\cdot du_{\sss F,E'},
\end{equation}
which was employed by Alexander to give a direct proof 
of the Schl\"{a}fli differential formula in the Euclidean case \cite{RA}.
Plug (\ref{equation_differential_theta_F}) into (\ref{equation_schlafli}), then
\begin{equation}
\label{equation_schlafli_E_F}
\kappa\cdot dV_k(P)=\frac{-1}{k-1}\sum_{F\subset E\subset P}
V_{k-2}(F) u_{\sss E,P}\cdot du_{\sss F,E},
\end{equation}
which will be useful in the proof of the following theorem. 

\begin{theorem}
\emph{(Main Theorem 3)}
\label{theorem_stress_volume}
Let $K(t)$ be a family of cell complexes in $M^d$
depending smoothly on a parameter $t$ and $K(0)=K$,
and $\omega$ be a $(k+1)$-stress $(k\geq 1)$ on $K$. Then
\[\sum_{G\in K^k}\omega(G)\,dV_k(G) = 0
\]
at $t=0$, where the sum is taken over all $k$-cells $G$ of $K$.
\end{theorem}

\begin{proof}
For $k=1$,
if $G$ is a $1$-cell of $K$ and $B$ is a vertex of $G$, 
let $u_{\sss B,G}$ be the inward unit normal to $G$ at $B$.
Then we have
$dV_1(G) = -\sum_{B\in G} u_{\sss B,G}\cdot dB$.
Taking the sum over all $1$-cells $G$ of $K$, we have
\begin{align*}
\sum_{G\in K^1}\omega(G)\,dV_1(G)
&=-\sum_{G\in K^1}\omega(G)\sum_{B\in G}u_{\sss B,G}\cdot dB        \\
&=-\sum_{\{B\}\in K^0}(\sum_{G\in K^1,B\in G}\omega(G)u_{\sss B,G})\cdot dB.
\end{align*}
As $\omega$ is a $2$-stress on $K$, by Definition \ref{defintion_k_stress}
$\sum_{G\in K^1,B\in G}\omega(G)u_{\sss B,G}$ is $0$ for each vertex $B$,
so the above formula is $0$ at $t=0$.

For $k\geq 2$ and $\omega$ is a $(k+1)$-stress on $K$ in $M^d$, 
we first consider the non-Euclidean case. Applying (\ref{equation_schlafli_E_F})
on each $G\in K^k$,
\begin{align*}
&\kappa\cdot(k-1)\sum_{G\in K^k}\omega(G)\,dV_k(G)        \\
&\quad= -\sum_{G\in K^k}\omega(G)\sum_{F\subset E\subset G}V_{k-2}(F)
u_{\sss E,G}\cdot du_{\sss F,E}               \\
&\quad= -\sum_{\{F,E|F\subset E\}}V_{k-2}(F)
(\sum_{G\in K^k,E\subset G}\omega(G)u_{\sss E,G})\cdot du_{\sss F,E}.
\end{align*}
As $\omega$ is a $(k+1)$-stress on $K$, by Definition \ref{defintion_k_stress}
$\sum_{G\in K^k,E\subset G}\omega(G)u_{\sss E,G}=0$ for each $(k-1)$-cell $E$ of $K$, 
so the above formula is $0$ at $t=0$. Since $\kappa\ne 0$ for the non-Euclidean
case, therefore $\sum_{G\in K^k}\omega(G)\,dV_k(G)=0$ at $t=0$.

For the Euclidean case, we show that the $(k+1)$-stress $\omega$
can be treated as a ``limit'' of some spherical $(k+1)$-stresses.
For any $r>0$,
embed $\mathbb{R}^d$ into $\mathbb{R}^{d+1}$, $x\hookrightarrow (x,r)$,
and let $\mathbb{S}^d_r$ be a $d$-dimensional sphere in $\mathbb{R}^{d+1}$ 
with radius $r$ and 
centered at the origin $O$ of $\mathbb{R}^{d+1}$.
By a radial projection (from the center of $\mathbb{S}^d_r$)
of $\mathbb{R}^d$ onto $\mathbb{S}^d_r$,
we obtain a family (with respect to $t$)
of spherical cell complexes $K_r(t)$ in $\mathbb{S}^d_r$ from
the Euclidean cell complexes $K(t)$ in $\mathbb{R}^d$. Denote $K_r(0)$ by $K_r$.
For each $k$-cell $G$ of $K$, denote by $G_r$ the corresponding
spherical cell of $K_r$,
and by $v_{\sss G,r}$ the altitude vector for the point $O$ 
with respect to the affine span of $G$.
We define a real-valued function $\omega_r$ on all the $k$-cells $G_r$ of $K_r$ by
\begin{equation}
\label{equation_radial_projection}
\omega_r(G_r):=\omega(G)\cdot\frac{\norm{v_{\sss G,r}}}{r}.
\end{equation}
By Definition \ref{defintion_k_stress} 
$\sum_{G\in K^k,E\subset G}\omega(G)u_{\sss E,G}=0$ for each $(k-1)$-cell $E$ of $K$.
Let $u'_{\sss E,G}$ be the orthogonal component of $u_{\sss E,G}$ 
that is perpendicular to the linear span of $E$ under the embedding,
then $\sum_{G\in K^k,E\subset G}\omega(G)u'_{\sss E,G}=0$ 
for each $(k-1)$-cell $E$ of $K$. 

It is easy to see that 
\[u'_{\sss E,G}=\frac{\norm{v_{\sss G,r}}}{\norm{v_{\sss E,r}}}\cdot u_{\sss E_r,G_r},
\]
so 
\[\sum_{G\in K^k,E\subset G}\omega(G)\norm{v_{\sss G,r}}\cdot u_{\sss E_r,G_r}=0,
\]
thus by (\ref{equation_radial_projection}) we have
$\sum_{G\in K^k,E\subset G}\omega_r(G_r)u_{\sss E_r,G_r}=0$ 
for each $(k-1)$-cell $E$ of $K$.
By Definition \ref{defintion_k_stress} 
$\omega_r$ is a $(k+1)$-stress
on $K_r$. Therefore, for any fixed $r$, we have
$\sum_{G_r\in K_r^k}\omega_r(G_r)\,dV_k(G_r)=0$ at $t=0$. 
As when $r\rightarrow\infty$, we have $\frac{\norm{v_{\sss G,r}}}{r}\rightarrow 1$,
so by (\ref{equation_radial_projection}) $\omega_r(G_r)$ converges to $\omega(G)$,
and $K_r(t)$ converges uniformly to $K(t)$ with respect to
small $t\geq 0$.
Thus $\sum_{G\in K^k}\omega(G)\,dV_k(G)=0$ at $t=0$.
This completes the proof.
\end{proof}

Particularly in the Euclidean case, but not in the non-Euclidean case,
Theorem \ref{theorem_stress_volume} implies the following property. 
Let $\omega$ be a $(k+1)$-stress $(k\geq 1)$ on a cell complex $K$ in $\mathbb{R}^d$,
then $K$ can be proportionally scaled with a factor $t$, with the \emph{same} 
$\omega$ as a $(k+1)$-stress. As the volumes of all $k$-faces of $K$ are scaled 
with a factor $t^k$, then by Theorem \ref{theorem_stress_volume} 
and taking the derivative at $t=1$, we have
\[\sum_{G\in K^k}\omega(G)V_k(G)=0.
\]

\begin{remark}
While Theorem \ref{theorem_stress_volume} was initially developed 
as a tool to prove Theorem \ref{theorem_mainone_r_s_h_n},
it is a much more general result.
Theorem \ref{theorem_stress_volume} establishes a correspondence
between the signs of tension-compression constraints
and the signs of $(k+1)$-stresses on $k$-faces of a cell complex,
which generalizes the well established correspondence of $k=1$.
Namely, if the signs of volume constraints agree with
the signs of \emph{any} $(k+1)$-stress on $k$-faces,
then the volumes of all $k$-faces are \emph{instantaneously} preserved at $t=0$.
If $(k+1)$-stresses can also be assigned in a continuous manner
over $t$ on the family of cell complexes,
then the volumes of all $k$-faces are preserved for small $t\ge 0$.
So to some extent, it justifies the \emph{physical} meaning of
$(k+1)$-stresses, which was first introduced
more of a mathematical concept for $k>1$ by Lee.
\end{remark}

It can be summarized as follows.

\begin{corollary}
\label{corollary_stress_volume}
Let $K(t)$ be a family of cell complexes in $M^d$
depending smoothly on a parameter $t$ and $K(0)=K$,
and $\omega_t$ be $(k+1)$-stresses on $K(t)$ in a continuous manner 
over $t\ge 0$ and $\omega_0=\omega$.
Then $K(t)$ cannot be a non-trivial deformation for small $t\ge 0$
under which the volumes of $k$-faces with negative signs of $\omega$ 
only decrease (resp. increase), and the volumes of $k$-faces 
with positive signs of $\omega$ only increase (resp. decrease).
Here by non-trivial it means that the volume of at least one $k$-face is non-constant.
\end{corollary}

\subsection{Proof of Theorem \ref{theorem_mainone_r_s_h_n}}

We begin with some basic notions.
Let $\Lambda(\mathbb{R}^{d+1})$
be the exterior algebra of $\mathbb{R}^{d+1}$.
An inner product on 
$\Lambda^{k}(\mathbb{R}^{d+1})$,
induced by the standard inner product on $\mathbb{R}^{d+1}$,
can be \emph{well defined} by
\begin{equation}
\label{equation_exterior_inner_product}
(r_{1}\wedge\cdots\wedge r_{k})\cdot(s_{1}\wedge\cdots\wedge s_{k})
:= \det(r_{i}\cdot s_{j})_{1\leq i,j\leq k},
\end{equation}
with extension by bilinearity,
where $r_i$ and $s_i$ are any $2k$ elements in $\mathbb{R}^{d+1}$,

The notions above can be extended to $\Lambda(\mathbb{R}^{d,1})$ in parallel,
with the exception that the inner product on $\Lambda^k(\mathbb{R}^{d,1})$
is not positive definite, but it is not a concern of this paper.
Particularly if $F$ is a $k$-simplex in $\mathbb{S}^d$ or $\mathbb{H}^d$ (recall that they
are embedded in $\mathbb{R}^{d+1}$ and $\mathbb{R}^{d,1}$ respectively)
and $B_1$, \dots, $B_{k+1}$
are the vertices, for convenience we introduce a new notation
\[\norm{F}:=|\det(B_i\cdot B_j)_{1\leq i,j\leq k+1}|^{1/2}.
\]
For the spherical case $\norm{F}$ is
simply $(k+1)!$ times the volume of the Euclidean $(k+1)$-simplex
whose vertices are $O$, $B_1$, \dots, $B_{k+1}$;
for the hyperbolic case, pseudo-volume.
With the volume interpretation of $\norm{F}$ in mind,
it will be very helpful for understanding the calculations involving $\norm{F}$
for the rest of this paper.

With the new notation, we have the following definition 
for a more general $(\mathbf{A},\alpha)$.

\begin{definition}
\label{definition_k_stress_alpha}
Let $(\mathbf{A},\alpha)$ be as in (\ref{equation_alpha_r_s_h})
where $\alpha$ is a 1-stress on $\mathbf{A}$,
but $\mathbf{A}$ is more general and may contain $m\geq n+2$ points
in general position in $M^n$.
For a given $k$ $(1\leq k\leq n)$ and each
simplicial $k$-face $F$ of $\mathbf{A}$,
define a $(k+1)$-stress $\alpha^{k+1}$
by $\alpha^{k+1}(F):=(\prod_{A_s\in F}\alpha_s)k!V_k(F)$
for the Euclidean case, and
$\alpha^{k+1}(F):=(\prod_{A_s\in F}\alpha_s)\norm{F}$
for the non-Euclidean case.
\end{definition}

\begin{remark}
Recall the discussion in Section \ref{section_stresses} about 
the relationship between our notion of $(k+1)$-stresses and
Lee's affine and linear $(k+1)$-stresses,
it is not hard to see the following general fact:
For a simplicial $k$-face $F$ of $\mathbf{A}$,
$\prod_{A_s\in F}\alpha_s$ corresponds 
to Lee's affine $(k+1)$-stress in the Euclidean case,
or to Lee's linear $(k+1)$-stress in the non-Euclidean case.
For notational reasons, we use $\alpha^{k+1}$ to denote the $(k+1)$-stress
obtained by multiplying $\alpha$ with itself for $k+1$ times
and then \emph{normalized} by a volume factor,
rather than taking the value of $\prod_{A_s\in F}\alpha_s$ directly.
\end{remark}

Then by Theorem \ref{theorem_stress_volume},
we immediately have the following fact.

\begin{proposition}
\label{proposition_stress_volume_A}
Let $(\mathbf{A},\alpha)$ be as in (\ref{equation_alpha_r_s_h})
where $\alpha$ is a 1-stress on $\mathbf{A}$,
and $\alpha^{k+1}$ be a $(k+1)$-stress on $\mathbf{A}$
as in Definition \ref{definition_k_stress_alpha}.
Then by Theorem \ref{theorem_stress_volume}
\begin{equation}
\label{equation_stress_volume_A}
\sum_{F\subset \mathbf{A}, \dim(F)=k}
\alpha^{k+1}(F)\,dV_k(F) = 0
\end{equation}
holds at $t=0$.
\end{proposition}

Then it leads to the proof of Theorem \ref{theorem_mainone_r_s_h_n}.

\begin{proof} [Proof of Theorem \ref{theorem_mainone_r_s_h_n}]
For a given $k$, $1\leq k\leq n$,
let $\alpha^{k+1}$ be the $(k+1)$-stress on $\mathbf{A}$ as above in
Proposition \ref{proposition_stress_volume_A}.
As $\mathbf{A}(t)$ is confined in $M^n$, so
$\mathbf{A}(t)$ is degenerate for $t\geq 0$.
This allows us to assign 1-stresses $\alpha_t$ on $\mathbf{A}(t)$
in a continuous manner over $t$.
Therefore by Definition \ref{definition_k_stress_alpha}
we can assign $(k+1)$-stresses $\alpha^{k+1}_t$ on $\mathbf{A}(t)$
continuously over $t$ as well.
As the signs of volume constraints of $G_{n,k}$ and $G'_{n,k}$ agree
with the signs of $\alpha^{k+1}_t$ (including the opposite of) for small $t\ge 0$,
Theorem \ref{theorem_mainone_r_s_h_n} is just a special case
of Corollary \ref{corollary_stress_volume}.
This completes the proof.
\end{proof}

To see if Theorem \ref{theorem_mainone_r_s_h_n} can be improved to claim
that $G_{n,k}$ and $G'_{n,k}$
are rigid in $M^n$,
check Remark \ref{remark_theorem_mainone}.

\subsection{A key definition $g_{\sss F}(P,Q)$}
\label{section_key_definition}

As Question \ref{question_volume_preserve} is settled for case $d=n$
primarily using a new
property of $k$-stresses (Theorem \ref{theorem_stress_volume}),
we plan to apply similar techniques
for the more general case $d\geq n+1$.
However, unlike the case $d=n$, 
for $d\geq n+1$ when $\mathbf{A}(t)$ moves in $M^d$ and
is not confined in a lower dimensional $M^n$, 
there is no $1$-stress on $\mathbf{A}(t)$,
and therefore no $(k+1)$-stress on $\mathbf{A}(t)$ for $t>0$.
By applying Theorem \ref{theorem_stress_volume},
though we can still show that for both
$G_{n,k}$ and $G'_{n,k}$ the volumes of all $k$-faces are
\emph{instantaneously} preserved at $t=0$, it is a weaker result
than what we are looking for, i.e., like Theorem \ref{theorem_mainone_r_s_h_n_plus_1}.
Thus Theorem \ref{theorem_stress_volume}
alone is not enough for our purposes.

To fix this issue, in Section \ref{section_invariant}, 
for each $k$-stress $\omega$ on a cell complex in $M^d$,
we discover an invariant $c_k(\omega)$ associated with $\omega$.
This is one of the most important results of this paper,
and this invariant leads to both the formulation and proof 
of Theorem \ref{theorem_mainone_r_s_h_n_plus_1}.
In this section we first introduce a notion $g_{\sss F}(P,Q)$
in Definition \ref{definition_partial_derivative_g_theta},
an important step for introducing the invariant $c_k(\omega)$.
We also address the properties of $g_{\sss F}(P,Q)$ in detail,
which is of interest by its own right.

Consider a $k$-dimensional simplex $F$ and two points $P$ and $Q$ in $M^d$,
and denote by $\hat{F}$ the $(k+2)$-dimensional simplex in $M^d$ which is
the \emph{join} of $F$ with the segment $PQ$.
Also let $\theta_F$ be the dihedral angle of $\hat{F}$ at face $F$.
Assume $\hat{F}$ is non-degenerate,
then all edge lengths of $\hat{F}$ can vary independently of each other,
thus $\theta_F$ can vary in
such a manner that the distances between any pair of vertices of $\hat{F}$
are preserved except between $P$ and $Q$.
It follows that $V_{k+2}(\hat{F})$ can be treated as a function of a 
\emph{single} variable $\theta_F$, and we write the differential as
$dV_{k+2}(\hat{F})/d\theta_F$.%
\footnote{It should not be confused with another similar notion
that treats all the dihedral angles of $\hat{F}$ as independent variables
in the non-Euclidean case. 
}

For the non-Euclidean case,
let $P'$ (resp. $Q'$) be the vertical projection of point $P$ (resp. $Q$)
on the linear span of $F$.
Then $(P-P')\cdot (Q'-P')=0$ and $(Q-Q')\cdot (Q'-P')=0$.
So if $\theta_F$ varies while all edge lengths of $\hat{F}$ are fixed 
except between $P$ and $Q$, then
\begin{align*}
d\overrightarrow{PQ}^2
&=d((Q-Q')+(Q'-P')-(P-P'))^2=-2d((P-P')\cdot (Q-Q'))        \\
&=-2\norm{P-P'}\cdot\norm{Q-Q'}d\cos\theta_F
=2\norm{P-P'}\cdot\norm{Q-Q'}\cdot\sin\theta_F d\theta_F      \\
&=2\cdot\frac{\norm{\hat{F}}}{\norm{F}}d\theta_F,
\end{align*}
where the second step is because the squared terms are constants when $\theta_F$ varies,
and the last step uses the volume interpretation of $\norm{F}$ and $\norm{\hat{F}}$.
Therefore we obtain
\begin{equation}
\label{equation_partial_derivative_theta_edge}
\frac{dV_{k+2}(\hat{F})}{d\theta_F} =
2\cdot\frac{\norm{\hat{F}}}{\norm{F}}
\cdot\partial_{\overrightarrow{PQ}^2}V_{k+2}(\hat{F}),
\end{equation}
where $\partial_{\overrightarrow{PQ}^2}$ is the partial derivative
with respect to $\overrightarrow{PQ}^2$ with all other edge lengths
of $\hat{F}$ fixed.

This interpretation of $dV_{k+2}(\hat{F})/d\theta_F$
can be easily extended to 
$k$-dimensional convex polytope $F$ that is not necessarily simplicial. 
Consider two points $P$ and $Q$
in $M^d$ such that the segment $PQ$ is in general position with respect to $F$,
denote by $\hat{F}$ the $(k+2)$-dimensional polytope in $M^d$
which is the join of $F$ with the segment $PQ$,
and by $\theta_F$ the dihedral angle of $\hat{F}$ at face $F$.
Here it is not crucial for $\hat{F}$ to be a convex polytope in the strict sense,
and some degeneracy is allowed as long as
$V_{k+2}(\hat{F})$ and $\theta_F$ can be properly defined.
Same as the simplicial case above,
$V_{k+2}(\hat{F})$ can be treated as a function of a 
single variable $\theta_F$.

Now we give a key definition, a new definition introduced in this paper.

\begin{definition}
\label{definition_partial_derivative_g_theta}
Let $F$ be a $k$-dimensional convex polytope in $M^d$
and $\hat{F}$, $\theta_F$
be as above. If $\theta_F$ varies while all edge lengths of $\hat{F}$ are fixed
except between $P$ and $Q$,
then define $g_{\sss F}: M^d\times M^d \rightarrow \mathbb{R}$ by
\begin{equation}
\label{equation_partial_derivative_g_theta}
g_{\sss F}(P,Q):= (k+2)!\,\frac{d V_{k+2}(\hat{F})}{d\theta_F}.
\end{equation}
Also set $g_{\sss\varnothing}(P,Q)=1$.
\end{definition}

For a $k$-polytope $F$,
note that if we decompose it into simplices $F_1,\dots,F_m$,
then it induces a decomposition of $\hat{F}$ into
$\hat{F_1},\dots,\hat{F_m}$. As $\theta_{F_i}=\theta_F$ for
each $F_i$, so by (\ref{equation_partial_derivative_g_theta})
we immediately have the following fact.

\begin{lemma}
\label{lemma_g_decompose}
If $F$ is decomposed into simplices $F_1,\dots,F_m$, then
$g_{\sss F}=\sum_i g_{\sss F_i}$.
\end{lemma}

For a $k$-dimensional simplex $F$ in $\mathbb{S}^d$ or $\mathbb{H}^d$,
by (\ref{equation_partial_derivative_g_theta}) and 
(\ref{equation_partial_derivative_theta_edge}) we have
\begin{equation}
\label{equation_partial_derivative_g_edge}
g_{\sss F}(P,Q) =
2\cdot(k+2)!\frac{\norm{\hat{F}}}{\norm{F}}
\cdot\partial_{\overrightarrow{PQ}^2}V_{k+2}(\hat{F}).
\end{equation}
To give a explicit formula for $g_{\sss F}(P,Q)$,
we introduce the following notation.
Let $G$ be a $k$-simplex in $\mathbb{S}^d$ or $\mathbb{H}^d$ 
and $B_1$, \dots, $B_{k+1}$ be the vertices, then define
\begin{equation}
\label{equation_regress_coeff}
R^i_j(G)
:=(-1)^{i+j}\frac{(B_1\wedge\cdots\wedge\hat{B_i}\wedge\cdots\wedge B_{k+1})
\cdot(B_1\wedge\cdots\wedge\hat{B_j}\wedge\cdots\wedge B_{k+1})}
{(B_1\wedge\cdots\wedge\hat{B_i}\wedge\cdots\wedge B_{k+1})^2}.
\end{equation}
Roughly speaking,
for $i\neq j$, if $B_i$ is projected onto 
the linear span of $B_1$, \dots, $\hat{B_i}$, \dots, $B_{k+1}$,
and expressed as $\sum_{s\ne i}\beta_s B_s$,
then $-\beta_j$ is $R^i_j(G)$;
and $R^i_i(G)=1$.
For notational reasons,
if $E_j$ is the $(k-1)$-face $G\setminus\{B_j\}$,
for $i\ne j$ we also define $R^i_j(E_j):=0$.

We now give the explicit formula for $g_{\sss F}(P,Q)$ 
in $\mathbb{S}^d$ or $\mathbb{H}^d$ when $F$ is simplicial.

\begin{lemma}
\label{lemma_prod_s_h}
Let $G$ be a $k$-simplex in $\mathbb{S}^d$ or $\mathbb{H}^d$
of constant curvature $\kappa$
and $B_1$, \ldots, $B_{k+1}$ be the vertices.
Also let $E_i$ be the $(k-1)$-face $G\setminus\{B_i\}$,
and $F_{ij}$ be the $(k-2)$-face $G\setminus\{B_i,B_j\}$.
Then for $k\geq 2$
\[ \kappa\cdot \norm{F_{12}}\,g_{\sss F_{12}}(B_1,B_2)
=k(k-2)!\sum_{i<j}R^{ij}_{12}(G)\norm{F_{ij}}V_{k-2}(F_{ij}),
\]
where for $i\ne j$ and $s\ne t$,
\begin{equation}
\label{equation_coefficient_double_index}
R^{ij}_{st}(G):=R^i_s(G)R^j_t(E_i)+R^j_s(G)R^i_t(E_j).
\end{equation}
Particularly, $R^{st}_{st}(G)=1$, and $R^{it}_{st}(G)=R^{i}_{s}(G)$.
\end{lemma}

We want to point out that
as $g_{\sss F_{st}}$ is symmetric on $B_s$ and $B_t$,
we have $R^{ij}_{st}(G)=R^{ij}_{ts}(G)$,
although it is not so obvious to see from
(\ref{equation_coefficient_double_index}) itself.

While we defer the proof of Lemma \ref{lemma_prod_s_h}
to Section \ref{section_proof_of_lemma_prod_s_h},
we give a direct consequence of Lemma \ref{lemma_prod_s_h} here,
a new version
of Schl\"{a}fli differential formula on simplices based on edge lengths.

\begin{proposition}
\label{proposition_schlafli_edge}
\emph{(Schl\"{a}fli differential formula on simplices based on edge lengths)}
Let $G$ be a $k$-simplex in $\mathbb{S}^d$ or $\mathbb{H}^d$
and $B_1$, \dots, $B_{k+1}$ be the vertices,
and $F_{ij}$ be the $(k-2)$-face $G\setminus \{B_i,B_j\}$.
Then
\[2\cdot k!\,\norm{G}\,dV_k(G)
=\sum_{i<j}\norm{F_{ij}}\,g_{\sss F_{ij}}(B_i,B_j)
\,d\overrightarrow{B_iB_j}^2,
\]
where the explicit formula of
$g_{\sss F_{ij}}(B_i,B_j)$
is given in Lemma \ref{lemma_prod_s_h}.
\end{proposition}

\begin{proof} 
Apply (\ref{equation_partial_derivative_g_edge}) and chain rule.
\end{proof}

Particularly for $g_{\sss F}$ when $F$ is a single point $B$, we have the following.

\begin{corollary}
\label{corollary_g_B}
Let $B$, $P$, $Q$ be three points in $\mathbb{S}^d$ or $\mathbb{H}^d$ 
of constant curvature $\kappa$, then
\begin{equation} 
\label{equation_g_B}
g_{\sss B}(P,Q)=\frac{2}{1+\kappa P\cdot Q}
\overrightarrow{PB}\cdot\overrightarrow{QB}.
\end{equation}
\end{corollary}

\begin{proof}
Let $G$ be a 2-simplex and $B_1=P$, $B_2=Q$, $B_3=B$ be the vertices.
Then by Lemma \ref{lemma_prod_s_h} we have
\[\kappa\cdot g_{\sss B_3}(B_1,B_2)=2(1+R^3_1(G)+R^3_2(G)).
\]
Multiplying $(B_1\wedge B_2)^2$ on both sides, 
and applying (\ref{equation_regress_coeff}) 
and (\ref{equation_exterior_inner_product}),
we have
\begin{align*}
&\kappa\cdot g_{\sss B_3}(B_1,B_2)(B_1\wedge B_2)^2   \\
&\quad=2((B_1\wedge B_2)^2+(B_1\wedge B_2)\cdot(B_2\wedge B_3)
-(B_1\wedge B_2)\cdot(B_1\wedge B_3))      \\
&\quad=2((B_1\wedge B_2)^2
+(B_1\cdot B_2)(B_2\cdot B_3)-B_2^2(B_1\cdot B_3)     \\
&\quad\quad-B_1^2(B_2\cdot B_3)+(B_1\cdot B_2)(B_1\cdot B_3))     \\
&\quad=2((B_1\wedge B_2)^2
-(\frac{1}{\kappa}-B_1\cdot B_2)(B_1\cdot B_3+B_2\cdot B_3)).
\end{align*}
As $(B_1\wedge B_2)^2
=B_1^2B_2^2-(B_1\cdot B_2)^2
=(\frac{1}{\kappa}-B_1\cdot B_2)(\frac{1}{\kappa}+B_1\cdot B_2)$,
factor out $(\frac{1}{\kappa}-B_1\cdot B_2)$ from above we have
\[(1+\kappa B_1\cdot B_2) g_{\sss B_3}(B_1,B_2)
=2(\frac{1}{\kappa}+B_1\cdot B_2-B_1\cdot B_3-B_2\cdot B_3)
=2\overrightarrow{B_1B_3}\cdot\overrightarrow{B_2B_3},
\]
which finishes the proof.
\end{proof}

Note that as $P,Q\rightarrow B$ and $\kappa P\cdot Q\rightarrow 1$,
$g_{\sss B}(P,Q)\sim\overrightarrow{PB}\cdot\overrightarrow{QB}$ 
is approximately
the Riemannian metric at $B$, with the difference that
$g_{\sss B}$ is defined on the whole $\mathbb{S}^d$ or $\mathbb{H}^d$
instead of on the tangent space at $B$.
The positive definiteness of $g_{\sss B}$ will be addressed
in Section \ref{section_positive_kernel}.

If $B$, $P$, $Q$ are three points in $\mathbb{R}^d$, 
by Definition \ref{definition_partial_derivative_g_theta} we have
\begin{equation}
\label{equation_g_B_Euclidean}
g_{\sss B}(P,Q)
=\norm{\overrightarrow{PB}}\cdot\norm{\overrightarrow{QB}}\cdot\frac{d\sin\theta_B}{d\theta_B}
=\norm{\overrightarrow{PB}}\cdot\norm{\overrightarrow{QB}}\cdot\cos\theta_B
=\overrightarrow{PB}\cdot\overrightarrow{QB},
\end{equation}
which can also be viewed as a limit of the spherical case (\ref{equation_g_B}).

%

\subsection{An invariant of $k$-stress}
\label{section_invariant}

Now we are ready to state a key result of this paper,
which leads to both the formulation and proof 
of Theorem \ref{theorem_mainone_r_s_h_n_plus_1}.

\begin{theorem}
\emph{(Main Theorem 4)}
\label{theorem_invariant_r_s_h}
Let $K$ be a cell complex in $M^d$ of constant curvature $\kappa$
and $\omega$ be a $k$-stress on $(k-1)$-faces of $K$ for $k\ge 1$. 
Then as long as $g_{\sss F}(P,Q)$ is properly defined for each $F\in K^{k-1}$,
we have
\begin{equation}
\label{equation_invariant_r_s_h}
\sum_{F\in K^{k-1}}\omega(F)\,g_{\sss F}(P,Q) = c_k(\omega),
\end{equation}
where $c_k(\omega)$
is an invariant independent of the choice of points $P,Q\in M^d$.
And for the non-Euclidean case,
\begin{equation}
\label{equation_invariant_s_h}
c_k(\omega) =\kappa (k+1)(k-1)!
\sum_{F\in K^{k-1}}\omega(F)\,V_{k-1}(F).
\end{equation}
\end{theorem}

We first give some examples to illustrate Theorem \ref{theorem_invariant_r_s_h}.

\begin{example}
\label{example_invariant_n_plus_1}
Let $K$ be a cell complex in $\mathbb{S}^n$ whose top dimensional cells
form a decomposition (not necessarily simplicial) of $\mathbb{S}^n$. Then the canonical $(n+1)$-stress
$\omega$ on $K$ can be defined by $\omega(F)=1$ for all $n$-dimensional faces
$F$ of $K$.
Then by (\ref{equation_invariant_s_h}) we have
$c_{n+1}(\omega)=(n+2)n!V_n(\mathbb{S}^n)$,
which is always positive no matter how the vertices of $K$ are positioned.
\end{example}

\begin{example}
Let $\omega$ be a 1-stress on a finite points $x_0, \dots, x_m$ in $\mathbb{R}^1$,
and $\omega(x_i)>0$ for $i>0$ and $\omega(x_0)=-1$.
First, if we treat $x_i$ for $i>0$ as a value taken by a discrete \emph{random variable} $X$ with 
probability $\omega(x_i)$, then $x_0$ can be treated as the \emph{mean} of $X$.
And second, in (\ref{equation_invariant_r_s_h}) set $k=1$ and $P=Q=0$,
then by applying (\ref{equation_g_B_Euclidean}) we have
 $c_1(\omega)=\sum_{i>0}\omega(x_i)x_i^2-x_0^2$,
which is the \emph{variance} $var(X)$ of random variable $X$.
A continuous analogue can be easily generalized for continuous random variables.
\end{example}

Now we are ready to prove Theorem \ref{theorem_invariant_r_s_h}.

\begin{proof} [Proof of Theorem \ref{theorem_invariant_r_s_h}]
We first consider the case that $K$ is simplicial in $\mathbb{S}^d$ or $\mathbb{H}^d$,
which is the most important step of the proof.
For a $(k-1)$-face $F$ of $K$,
denote by $\hat{F}$ the $(k+1)$-dimensional simplex
which is the join of $F$ with segment $PQ$.
By Lemma \ref{lemma_prod_s_h}, we can view $g_{\sss F}(P,Q)$ 
as a weighted sum of the volumes of all $(k-1)$-faces $F'$ of $\hat{F}$,
and particularly, the weight on $F$ is independent of the choice of $P$ and $Q$.
There are four types of $F'$: 
$\{P,Q\}\subset F'$; $P\in F', Q\not\in F'$; $P\not\in F', Q\in F'$; and $F'=F$.
When summing over all $(k-1)$-faces $F$ of $K$
on the left side of (\ref{equation_invariant_r_s_h}),
applying Lemma \ref{lemma_prod_s_h}
and (\ref{equation_regress_coeff})
with some linear algebra,
it can be shown that
for any given $F'$ of the first three types,
the sum of the weights on $F'$ is $0$. So only the $4$-th type of terms
are left, and
(\ref{equation_invariant_s_h}) immediately follows.

For the more general case that $K$ is not necessarily simplicial 
in $\mathbb{S}^d$ or $\mathbb{H}^d$,
by a barycentric subdivision of all the cells of $K$ with dimension $k-1$ and lower,
and ignoring all the cells with dimension $k$ and above,
we obtain a simplicial cell complex $K'$ with dimension $k-1$.
Note that any $(k-1)$-simplicial face $F'$ of $K'$ is obtained from the
decomposition of a $(k-1)$-cell $F$ of $K$.
We define a real-valued function $\omega'$
on $F'$ of $K'$ by $\omega'(F'):=\omega(F)$.
It can be shown that $\omega'$ is a $k$-stress on $K'$
by using Definition \ref{defintion_k_stress} to verify
the following two types of $(k-2)$-simplicial faces $G'$ of $K'$:
the first type of $G'$ is 
part of a $(k-2)$-cell $G$ of $K$, so it automatically satisfies the condition 
in Definition \ref{defintion_k_stress};
the second type of $G'$ is introduced by the decomposition of a $(k-1)$-cell $F$ of $K$
but not on the boundary of $F$, so $G'$ is shared by exactly two
$(k-1)$-simplicial faces of $K'$ who have the opposite inward unit normals at
their common facet $G'$, and therefore satisfies the condition 
in Definition \ref{defintion_k_stress} as well.

Applying the facts
(1) $\omega'$ is a $k$-stress on $K'$,
(2) the simplicial version
of (\ref{equation_invariant_r_s_h}) and (\ref{equation_invariant_s_h})
we just proved above,
and (3) the formula $g_{\sss F}=\sum_i g_{\sss F_i}$ from Lemma \ref{lemma_g_decompose},
we prove that
(\ref{equation_invariant_r_s_h}) and (\ref{equation_invariant_s_h})
still hold for the case that when $K$ is a not necessarily simplicial
in $\mathbb{S}^d$ or $\mathbb{H}^d$.

Finally, for the Euclidean case, let $\omega$ be a $k$-stress on $K$.
For any $r>0$,
embed $\mathbb{R}^d$ into $\mathbb{R}^{d+1}$, $x\hookrightarrow (x,r)$,
and let $\mathbb{S}^d_r$ be a $d$-dimensional sphere in $\mathbb{R}^{d+1}$ 
with radius $r$ and 
centered at the origin $O$ of $\mathbb{R}^{d+1}$.
By a radial projection (from the center of $\mathbb{S}^d_r$)
of $\mathbb{R}^d$ onto $\mathbb{S}^d_r$
(see (\ref{equation_radial_projection}) and the nearby discussion),
it induces a $k$-stress in $\mathbb{S}^d_r$.
Taking $r\rightarrow\infty$ in the spherical case,
we prove that $\sum_{F\in K^{k-1}}\omega(F)\,g_{\sss F}(P,Q)$
of (\ref{equation_invariant_r_s_h})
is independent of $P$ and $Q$
in the Euclidean case as well.
This completes the proof.
\end{proof}

Particularly for a general $(\mathbf{A},\alpha)$, we have the following definition.

\begin{definition}
\label{defintion_invariant_c_k}
Let $(\mathbf{A},\alpha)$ be as in (\ref{equation_alpha_r_s_h})
where $\alpha$ is a 1-stress on $\mathbf{A}$,
but $\mathbf{A}$ is more general and may contain $m\geq n+2$ points
in general position in $M^n$.
Also let $\alpha^k$ be a $k$-stress on $\mathbf{A}$
as in Definition \ref{definition_k_stress_alpha}.
Then by Theorem \ref{theorem_invariant_r_s_h} we define 
a sequence of invariants
$c_1(\alpha^1), \dots, c_{n+1}(\alpha^{n+1})$ for $(\mathbf{A},\alpha)$.
Also set $c_0(\alpha^0)=1$.
For the non-Euclidean case by Theorem \ref{theorem_invariant_r_s_h}
we have
\begin{equation}
\label{equation_invariant_s_h_A}
c_k(\alpha^k) =\kappa (k+1)(k-1)!
\sum_{F\subset \mathbf{A},\dim(F)=k-1}(\prod_{A_s\in F}\alpha_s)
\norm{F}
\,V_{k-1}(F).
\end{equation}
\end{definition}

\begin{remark}
\label{remark_invariant_n_plus_1}
When $\mathbf{A}$ contains exactly $m=n+2$ points,
for the non-Euclidean case, by (\ref{equation_invariant_s_h_A})
$c_{n+1}(\alpha^{n+1})$ vanishes unless
$\mathbf{A}$ is not confined
in any open half sphere in the spherical case
(see Example \ref{example_invariant_n_plus_1}
for a case that it does not vanish);
$c_{n+1}(\alpha^{n+1})$ also vanishes in the Euclidean case
as a limit of the spherical case.
Using the special case of $m=n+2$,
with the proof skipped,
the same conclusion can be proved for $m>n+2$
or even when $(\mathbf{A},\alpha)$ is distributed in a continuous manner.
\end{remark}

Note that if $\alpha_i>0$ for all $i$, which can only happen in the spherical case,
then all $c_k(\alpha^k)$ are positive, including $c_{n+1}(\alpha^{n+1})$.

\begin{corollary}
\label{corollary_invariant_r_s_h_A}
Let $\alpha^k$ be the $k$-stress
on $\mathbf{A}$ as above in Definition \ref{defintion_invariant_c_k}.
Then for $k\leq n$ and any $i\ne j$, we also have
\begin{equation}
\label{equation_invariant_r_s_h_A}
\sum_{\substack{F\subset\mathbf{A}\setminus\{A_i,A_j\}\\ \dim(F)=k-1}}
\alpha^k(F)\,g_{\sss F}(A_i,A_j)=c_k(\alpha^k).
\end{equation}
\end{corollary}

\begin{proof}
Following essentially the same proof of Theorem \ref{theorem_invariant_r_s_h},
the non-Euclidean case can be proved by applying
Lemma \ref{lemma_prod_s_h} and (\ref{equation_regress_coeff}),
and the Euclidean case can be proved by treating it as a limit of the spherical case.
\end{proof}

\subsection{Proof of Lemma \ref{lemma_prod_s_h}}
\label{section_proof_of_lemma_prod_s_h}

The proof of Lemma \ref{lemma_prod_s_h} is mainly computational,
the reader not interested in technicalities 
can skip this section for now without missing the flow of the paper.

Let $G$ be a $k$-simplex in $\mathbb{S}^d$ or $\mathbb{H}^d$
of constant curvature $\kappa$
and $B_1$, \ldots, $B_{k+1}$ be the vertices.
The main idea to prove Lemma \ref{lemma_prod_s_h}
is to compute $dV_k(G)$ in two different ways,
and to compare the coefficients of the
outcomes. One way is to expand $dV_k(G)$ as a linear sum of
$d\overrightarrow{B_iB_j}^2$ by using
(\ref{equation_partial_derivative_g_edge}),
and the other is to expand $dV_k(G)$ using
the Schl\"{a}fli differential formula (\ref{equation_schlafli}).

Let $E_i$ be the $(k-1)$-face $G\setminus\{B_i\}$ of $G$,
$F_{ij}$ the $(k-2)$-face that can be described as
an intersection $E_i\cap E_j$, and $\theta_{ij}$ the dihedral angle between
$E_i$ and $E_j$.
Also let $e_i$ be the inward unit normal to $G$
along the $(k-1)$-face $E_i$,
and $f_{ij}$ the inward unit normal to $E_i$
along the $(k-2)$-face $F_{ij}$. It is obvious that
the angle between $e_i$ and $f_{ij}$ is $\pi/2$, and therefore
$e_i\cdot f_{ij}=0$.
For $i\ne j$, recall (\ref{equation_differential_theta_F}) that
$d\theta_{ij}= -e_i\cdot df_{ij} - e_j\cdot df_{ji}$.
As $e_i\cdot f_{ij}=0$, we have
$e_i\cdot df_{ij} + f_{ij}\cdot de_i = 0$,
therefore
\begin{equation}
\label{equation_differential_theta}
d\theta_{ij}= f_{ij}\cdot de_i + f_{ji}\cdot de_j,
\end{equation}
which will be very useful in the proof of Lemma \ref{lemma_prod_s_h}.

Recall (\ref{equation_regress_coeff}) and the nearby interpretation of
$R^i_j(G)$, easy to see that
$\sum_s R^i_s(G)B_s$ is the altitude vector for point $B_i$ with respect to 
the linear span of $B_1$, \dots, $\hat{B_i}$, \dots, $B_{k+1}$.
As the norm of the altitude vector is $\frac{\norm{G}}{\norm{E_i}}$,
normalizing the vector we have
\begin{equation}
\label{equation_inward_normal_e}
e_i=\frac{\norm{E_i}}{\norm{G}}
\sum_s R^i_s(G)B_s.
\end{equation}
And similarly, because by definition for $i\ne j$ we have $R^j_i(E_i)=0$, thus
\begin{equation}
\label{equation_inward_normal_f}
f_{ij}=\frac{\norm{F_{ij}}}{\norm{E_i}}
\sum_s R^j_s(E_i)B_s.
\end{equation}

\begin{proof} [Proof of Lemma \ref{lemma_prod_s_h}]
For $k\geq 2$, by (\ref{equation_coefficient_double_index}),
Lemma \ref{lemma_prod_s_h} is equivalent to proving
\begin{equation}
\label{equation_partial_derivative_equivalent}
\kappa\cdot \norm{F_{12}}\,g_{\sss F_{12}}(B_1,B_2)
=k(k-2)!\sum_{i\ne j}R^i_1(G)R^j_2(E_i)
\norm{F_{ij}}V_{k-2}(F_{ij}).
\end{equation}

In the rest of the proof, we assume the vertices of $G$ are moving
in such a manner that $B_2$ is the only vertex that is moving,
and $\overrightarrow{B_2B_i}^2$ are preserved for $3\leq i\leq k+1$.
Under this assumption, $B_1\cdot dB_2$
is the only non-zero term of the form
$B_i\cdot dB_j$ for $1\leq i,j\leq k+1$.

So by
(\ref{equation_partial_derivative_g_edge}),
we have
\begin{equation}
\label{equation_differential_delta}
\begin{split}
&\kappa\cdot k!\,\norm{G}\,dV_k(G)      \\
&\quad=\frac{\kappa}{2}\cdot\norm{F_{12}}\,g_{\sss F_{12}}(B_1,B_2)
\,d\overrightarrow{B_1B_2}^2          \\
&\quad=-\kappa\cdot\norm{F_{12}}\,g_{\sss F_{12}}(B_1,B_2)\, d(B_1\cdot B_2)
=-\kappa\cdot\norm{F_{12}}\,g_{\sss F_{12}}(B_1,B_2)\,B_1\cdot dB_2.  \\
\end{split}
\end{equation}

On the other hand, applying the Schl\"{a}fli differential formula
(\ref{equation_schlafli}), we have
\allowdisplaybreaks{
\begin{align*}
&\kappa\cdot (k-1)\,\norm{G}\,dV_k(G)             \\
&\quad= \norm{G}\sum_{i<j}V_{k-2}(F_{ij})\,d\theta_{ij}
= \norm{G}\sum_{i\ne j}V_{k-2}(F_{ij})\,f_{ij}\cdot de_i  \\
&\quad= \norm{G}\sum_{i\ne j}V_{k-2}(F_{ij})
\frac{\norm{F_{ij}}}{\norm{E_i}}
\sum_t R^j_t(E_i)B_t\cdot de_i       \\
&\quad= -\norm{G}\sum_{i\ne j}V_{k-2}(F_{ij})
\frac{\norm{F_{ij}}}{\norm{E_i}}
\sum_t R^j_t(E_i)e_i\cdot dB_t     \\
&\quad= -\norm{G}\sum_{i\ne j}V_{k-2}(F_{ij})
\frac{\norm{F_{ij}}}{\norm{E_i}}
\sum_t R^j_t(E_i)
\frac{\norm{E_i}}{\norm{G}}
\sum_s R^i_s(G)B_s\cdot dB_t       \\
&\quad= -\sum_{i\ne j}V_{k-2}(F_{ij})\norm{F_{ij}}
R^i_1(G)R^j_2(E_i)
B_1\cdot dB_2,
\end{align*}
}
where the second step is by applying (\ref{equation_differential_theta});
the third step is by (\ref{equation_inward_normal_f});
the fourth step is because $e_i\cdot B_t=0$ when $i\ne t$
and $R^j_t(E_i)=0$ when $i=t$;
the fifth step is by (\ref{equation_inward_normal_e});
and the last step is because $B_1\cdot dB_2$
is the only non-zero term of the form
$B_i\cdot dB_j$ for $1\leq i,j\leq k+1$.
Comparing the last step above with
(\ref{equation_differential_delta}), we prove
(\ref{equation_partial_derivative_equivalent})
and complete the proof.
\end{proof}

\subsection{Proof of Theorem \ref{theorem_mainone_r_s_h_n_plus_1}}

The invariant $c_k(\alpha^k)$ (Definition \ref{defintion_invariant_c_k})
of $(\mathbf{A},\alpha)$ plays a role in both the formulation and proof of 
Theorem \ref{theorem_mainone_r_s_h_n_plus_1}.

Assume $\mathbf{A}(t)$ is in $M^d$ with $d\geq n+1$.
Let $A_0(t)$ in $M^d$ be the mirror reflection of $A_1(t)$ through
a lower dimensional $M^n$ that contains points $A_2(t),\ldots,A_{n+2}(t)$.
It is not hard to see that if $\mathbf{A}(t)$ varies smoothly over $t$,
then $A_0(t)$ varies smoothly as well.
We denote $\mathbf{A}(t)\cup \{A_0(t)\}$ by $\mathbf{A}^{\ast}(t)$.
By adding $A_0(t)$, we treat $\mathbf{A}^{\ast}(t)$
as a degenerate $(n+2)$-simplex in $M^d$.
So 
for each $i\ge 1$, $\alpha_i$ can be extended to a continuous
function $\alpha_i(t)$ with $\alpha_{i}(0)=\alpha_{i}$, such that
$\{\beta_0(t),\dots,\beta_{n+2}(t)\}$ is a 1-stress on $\mathbf{A}^{\ast}(t)$,
where $\beta_0(t)=\beta_1(t)=\frac{1}{2}\alpha_1(t)$
and $\beta_i(t)=\alpha_i(t)$ for $2\leq i\leq n+2$.

Denote $\{\alpha_1(t),\dots,\alpha_{n+2}(t)\}$ by $\alpha_t$
and $\{\beta_0(t),\dots,\beta_{n+2}(t)\}$ by $\beta_t$.
So Proposition \ref{proposition_stress_volume_A} can be
applied on $(\mathbf{A}^{\ast}(t),\beta_t)$ for $t\ge 0$.
Since $A_0(t)$ is the reflection of $A_1(t)$, so for each
$k$-face of $\mathbf{A}^{\ast}(t)$ that
contains point $A_1(t)$ but not $A_0(t)$, there is a
congruent
$k$-face that contains the same set of vertices except with
$A_1(t)$ replaced by $A_0(t)$.
These are the key ideas to prove the following result,
a final step before proving 
Theorem \ref{theorem_mainone_r_s_h_n_plus_1}.

\begin{proposition}
\label{proposition_inequality_r_s_h}
Let $\mathbf{A}(t)$, $\alpha_t$ and $(\mathbf{A}^{\ast}(t),\beta_t)$ be as above.
Assume $\mathbf{A}(t)$ varies smoothly for $t\ge 0$ in $M^d$ with $d\geq n+1$.
If $c_{k-1}(\alpha^{k-1})\ne 0$ and $A_0(t)\ne A_1(t)$ for small $t>0$, 
then for small $t>0$, for the non-Euclidean case
\begin{equation}
\label{equation_inequality_s_h}
2\cdot k! \sum_{\substack{G\subset \mathbf{A}(t)\\ \dim(G)=k}}
(\prod_{A_s(t)\in G}\alpha_s(t))
\,\norm{G}\,dV_k(G)
\sim -\frac{1}{4}\alpha_1^2 c_{k-1}(\alpha^{k-1})\,d\overrightarrow{A_0A_1}^2,
\end{equation}
and for the Euclidean case
\begin{equation}
\label{equation_inequality_r}
2\cdot (k!)^2\sum_{\substack{G\subset \mathbf{A}(t)\\ \dim(G)=k}}
(\prod_{A_s(t)\in G}\alpha_s(t))
\,V_k(G)\,dV_k(G)
\sim -\frac{1}{4}\alpha_1^2 c_{k-1}(\alpha^{k-1})\,d\overrightarrow{A_0A_1}^2.
\end{equation}
\end{proposition}

\begin{proof}
For the non-Euclidean case, for $(\mathbf{A}^{\ast}(t), \beta_t)$,
by Definition \ref{definition_k_stress_alpha}
let $\beta^{k+1}_t$ be the $(k+1)$-stress on
the $k$-faces $G$ of $\mathbf{A}^{\ast}(t)$ for $k\geq 0$
such that 
\[ \beta^{k+1}_t(G):=(\prod_{A_s(t)\in G}\beta_s(t))\norm{G}.
\]
Also set $\beta^0_t(\varnothing)=1$.
Now switch the index from $k+1$ to $k-1$ if $k\geq 1$,
by Theorem \ref{theorem_invariant_r_s_h},
$\beta^{k-1}_t$ has an associated invariant $c_{k-1}(\beta^{k-1}_t)$.
By applying (\ref{equation_invariant_s_h}), one sees that
as $t\rightarrow 0$, $c_{k-1}(\beta^{k-1}_t)$ converges to $c_{k-1}(\alpha^{k-1})$.
So if $c_{k-1}(\alpha^{k-1})\ne 0$, then for small $t>0$,
\allowdisplaybreaks{
\begin{align*}
&2\cdot k! \sum_{G\subset \mathbf{A}(t),\dim(G)=k}
(\prod_{A_s(t)\in G}\alpha_s(t))
\,\norm{G}\,dV_k(G)             \\
&\quad= 2\cdot k!
\sum_{\substack{G\subset \mathbf{A}^{\ast}(t),\dim(G)=k\\ \{A_0(t),A_1(t)\}\not\subset G}}
\beta^{k+1}_t(G)\,dV_k(G)             \\
&\quad= -2\cdot k!
\sum_{\substack{G\subset \mathbf{A}^{\ast}(t),\dim(G)=k\\ \{A_0(t),A_1(t)\}\subset G}}
\beta^{k+1}_t(G)\,dV_k(G)             \\
&\quad\sim -2\cdot k!
\sum_{\substack{G\subset \mathbf{A}^{\ast}(t),\dim(G)=k\\ \{A_0(t),A_1(t)\}\subset G}}
\beta^{k+1}_t(G)\,\partial_{\overrightarrow{A_0A_1}^2}V_k(G)\,d\overrightarrow{A_0A_1}^2     \\
&\quad= -\beta_0(t)\beta_1(t)
\sum_{\substack{F\subset \mathbf{A}^{\ast}(t)\setminus\{A_0(t),A_1(t)\}\\ \dim(F)=k-2}}
\beta^{k-1}_t(F)\,g_{\sss F}(A_0(t),A_1(t))
\,d\overrightarrow{A_0A_1}^2             \\
&\quad=-\beta_0(t)\beta_1(t)c_{k-1}(\beta^{k-1}_t)
\,d\overrightarrow{A_0A_1}^2         \\
&\quad\sim -\frac{1}{4}\alpha_1^2 c_{k-1}(\alpha^{k-1})
\,d\overrightarrow{A_0A_1}^2,
\end{align*}
}
where the second step is by
applying Proposition \ref{proposition_stress_volume_A} on $\mathbf{A}^{\ast}(t)$;
the third step is because 
$\norm{G}dV_k(G)\sim \norm{G}\partial_{\overrightarrow{A_0A_1}^2}V_k(G)\,d\overrightarrow{A_0A_1}^2$
for small $t>0$ when $\{A_0(t),A_1(t)\}\subset G$;
the fourth step is because (\ref{equation_partial_derivative_g_edge});
the fifth step is by applying Corollary \ref{corollary_invariant_r_s_h_A}
on $\mathbf{A}^{\ast}(t)$;
and the last step is because
$\beta_0(t)=\beta_1(t)=\frac{1}{2}\alpha_1(t)$ and
$c_{k-1}(\beta^{k-1}_t)$ converges to $c_{k-1}(\alpha^{k-1})$.
This completes the proof of (\ref{equation_inequality_s_h}).

For the Euclidean case, the only change we need to make,
is by Definition \ref{definition_k_stress_alpha}
let $\beta^{k+1}_t$ be the $(k+1)$-stress on
the $k$-faces $G$ of $\mathbf{A}^{\ast}(t)$
such that 
\[\beta^{k+1}_t(G):=(\prod_{A_s(t)\in G}\beta_s(t))k!V_k(G).
\]
Also set $\beta^0_t(\varnothing)=1$.
Following the same steps above,
we then prove (\ref{equation_inequality_r}).
\end{proof}

Note that the proof of Proposition \ref{proposition_inequality_r_s_h}
starts with a $(k+1)$-stress on $\mathbf{A}$,
but it is the invariant $c_{k-1}(\alpha^{k-1})$ of a $(k-1)$-stress, 
rather than the invariant $c_{k+1}(\alpha^{k+1})$, 
plays a role in both the formulation and proof.
The same applies to 
the following proof of Theorem \ref{theorem_mainone_r_s_h_n_plus_1},
which is a direct consequence of Proposition \ref{proposition_inequality_r_s_h}.

\begin{proof} [Proof of Theorem \ref{theorem_mainone_r_s_h_n_plus_1}]
We need only prove the non-Euclidean case,
as the Euclidean case can be proved similarly.
Also assume $c_{k-1}(\alpha^{k-1})>0$, 
as the case $c_{k-1}(\alpha^{k-1})<0$ is similar.

Let $A_0(t)$ and $\alpha_t=\{\alpha_1(t),\dots,\alpha_{n+2}(t)\}$ be the same as
in Proposition \ref{proposition_inequality_r_s_h}.
For $G_{n,k}$, now assume that the vertices are not always confined
in a lower dimensional $\mathbb{S}^n$ or $\mathbb{H}^n$ for small $t>0$,
then there exists arbitrarily small $t>0$ such that
$A_0(t)\ne A_1(t)$.%
\footnote{It is possible to still have infinitely many small $t>0$ such that
$A_0(t)=A_1(t)$, e.g., at the zeros of the function $e^{-1/t^2}\sin(1/t)$ near $t=0$, 
so we should be cautious about this kind of scenario. 
However, this is not a concern if $\mathbf{A}(t)$ is real analytic.
}
Therefore there exist arbitrarily small $t_1$ and $t_2$ with $0\leq t_1 < t_2$,
such that $A_0(t_1)=A_1(t_1)$,
but $A_0(t)\ne A_1(t)$ for any $t$ with $t_1<t<t_2$.
So $\alpha_{t_1}$ is a 1-stress on $\mathbf{A}(t_1)$, 
then $c_{k-1}(\alpha^{k-1}_{t_1})$ can be properly defined
by Definition \ref{defintion_invariant_c_k}.
As $c_{k-1}(\alpha^{k-1})>0$, 
so $t_1$ can be small enough such that 
$c_{k-1}(\alpha^{k-1}_{t_1})>0$ as well.

Then we can apply Proposition \ref{proposition_inequality_r_s_h}
to $\mathbf{A}(t)$ near $t=t_1$,
so for small  $t-t_1>0$,
\begin{equation}
\label{equation_inequaltiy_s_h_copy}
2\cdot k!
\sum_{\substack{F\subset \mathbf{A}(t)\\ \dim(F)=k}}(\prod_{A_s(t)\in F}\alpha_s(t))
\,\norm{F}\,dV_k(F)
\sim -\frac{1}{4}\alpha_1(t_1)^2\cdot c_{k-1}(\alpha^{k-1}_{t_1})\,d\overrightarrow{A_0A_1}^2.
\end{equation}
However, for each $k$-face $F$ of $G_{n,k}$,
by Definition \ref{definition_framework} we have
$(\prod_{A_s(t)\in F}\alpha_s)\,dV_k(F)\geq 0$ for $t\geq t_1\geq 0$,
which is a contradiction to (\ref{equation_inequaltiy_s_h_copy}).
So the vertices of $G_{n,k}$ must be confined in a lower dimensional
$\mathbb{S}^{n}$ or $\mathbb{H}^{n}$ for small $t\geq 0$.
Applying Theorem \ref{theorem_mainone_r_s_h_n},
we then show that $V_k(F)$
must be preserved for small $t\geq 0$.
This completes the proof.
\end{proof}

To see if Theorem \ref{theorem_mainone_r_s_h_n_plus_1}
can be improved to claim that $G_{n,k}$ or $G'_{n,k}$ is
rigid in $M^d$,
check Remark \ref{remark_theorem_mainone}.
More results and examples for cases $k=1$ and $2$ are given next.

\subsection{Tensegrity framework $G_{n,1}$}
\label{section_tensegrity_1}

In Theorem \ref{theorem_mainone_r_s_h_n_plus_1} for $k=1$,
it states that $G_{n,1}$ is rigid in $\mathbb{R}^d$ for $d\ge n+1$.
While this result is not new,
our theorem provides a new interpretation by using $c_0(\alpha^0)=1>0$.
Bezdek and Connelly \cite{BC} proved a stronger result that $G_{n,1}$ is
globally rigid in $\mathbb{R}^d$ for any $d\geq n$,
and with a little modification the spherical case
can be proved as well.
However, the methodology they used
cannot be directly applied to prove the hyperbolic case,
which mainly because the metric in $\mathbb{R}^{d,1}$ is not positive definite.
To our knowledge,
the following result we obtained in hyperbolic space
is new.

\begin{theorem}
\label{theorem_h_global_rigid}
$G_{n,1}$ is globally rigid in $\mathbb{H}^d$
for any $d\ge n$.
\end{theorem}

\begin{proof}
Recall (\ref{equation_alpha_r_s_h}) that $\sum_{i=1}^{n+2}\alpha_i A_i=0$.
Assume $\alpha_1$, \dots, $\alpha_m>0$ and
$\alpha_{m+1}$, \dots, $\alpha_{n+2}<0$.
Let $B_1$, \dots, $B_{n+2}$ be $n+2$ points in $\mathbb{H}^d$ that satisfy the
constraints of $G_{n,1}$. Namely, 
$\alpha_i\alpha_j\overrightarrow{B_iB_j}^2\geq
\alpha_i\alpha_j\overrightarrow{A_iA_j}^2$
for any $i\ne j$,
which is the same as
$\alpha_i\alpha_j B_i\cdot B_j\leq \alpha_i\alpha_j A_i\cdot A_j$.
Since $\alpha_1,\dots,\alpha_m > 0$, so there is a $f_1>0$ such that
$f_1\sum_{i=1}^{m}\alpha_i B_i$ is a point in $\mathbb{H}^d$; similarly,
there is a $f_2>0$ such that $-f_2\sum_{i=m+1}^{n+2}\alpha_i B_i$ is
a point in $\mathbb{H}^d$. Denote these two points by $D_1$ and $D_2$, and
let $\beta_i=f_1\alpha_i$ if $1\leq i\leq m$ and
$\beta_i=f_2\alpha_i$ if $m+1\leq i\leq n+2$.
As $\overrightarrow{D_2D_1}^2\geq 0$, so
\begin{align*}
0&\leq (D_1 - D_2)^2
=(\sum\beta_i B_i)^2 \leq(\sum\beta_i A_i)^2          \\
&=(f_1\sum_{i=1}^{m}\alpha_i A_i + f_2\sum_{i=m+1}^{n+2}\alpha_i
A_i)^2
=((f_1-f_2)\sum_{i=1}^m\alpha_i A_i)^2
\leq 0,
\end{align*}
where the third step is because
$\alpha_i\alpha_j B_i\cdot B_j\leq \alpha_i\alpha_j A_i\cdot A_j$
and $f_1,f_2>0$,
so $\beta_i\beta_j B_i\cdot B_j\leq \beta_i\beta_j A_i\cdot A_j$;
the fifth step is because $\sum_{i=1}^{n+2}\alpha_i A_i=0$;
the last step is because $\sum_{i=1}^m\alpha_i A_i$
is a multiple of a point in $\mathbb{H}^d$, so
$(\sum_{i=1}^m\alpha_i A_i)^2<0.$

Then $B_i\cdot B_j = A_i\cdot A_j$
holds for any $i\ne j$, and so $G_{n,1}$ is
globally rigid in $\mathbb{H}^d$.
\end{proof}

\subsection{2-tensegrity frameworks $G_{n,2}$ and $G'_{n,2}$}
\label{section_tensegrity_2}

In Theorem \ref{theorem_mainone_r_s_h_n_plus_1}
the sign of $c_1(\alpha^1)$ plays an important role in the case $k=2$.
In this section, we give a geometric interpretation of $c_1(\alpha^1)=0$,
which is amazingly simple as shown below.

\begin{proposition}
\label{proposition_common_sphere_hyperplane}
For the spherical (resp. hyperbolic) case,
$c_1(\alpha^1)=0$ if and only if
$A_1$, \dots, $A_{n+2}$ are affinely dependent
in $\mathbb{R}^{n+1}$ (resp. $\mathbb{R}^{n,1}$).
For the Euclidean case,
$c_1(\alpha^1)=0$ if and only if $A_1$, \dots, $A_{n+2}$ 
lie on a $(n-1)$-dimensional sphere in $\mathbb{R}^{n}$.
\end{proposition}

\begin{proof}
For the spherical (resp. hyperbolic) case,
by (\ref{equation_invariant_s_h_A})
we have $c_1(\alpha^1)=\kappa\cdot2\sum\alpha_{i}$.
Since $\sum\alpha_{i}A_{i}=0$,
so $c_1(\alpha^1)=0$ (the same as $\sum\alpha_{i}=0$) if and only if
$A_1$, \dots, $A_{n+2}$ are affinely dependent.

For the Euclidean case,
let $\mathbb{S}^{n-1}_1$ be a $(n-1)$-dimensional sphere in
$\mathbb{R}^{n}$ that contains points $A_{2}$, \dots, $A_{n+2}$;
$O_1$ be the center of the sphere and $r$ be the radius.
From (\ref{equation_g_B_Euclidean})
we have
$g_{\sss A_i}(P,Q)=\overrightarrow{PA_i}\cdot\overrightarrow{QA_i}$
in $\mathbb{R}^{n}$,
then by choosing $P=Q=O_1$ in (\ref{equation_invariant_r_s_h})
we have
$c_1(\alpha^1)=\sum\alpha_{i}\overrightarrow{O_1A_i}^2$.
Since $\sum\alpha_i=0$, so
$c_1(\alpha^1)=\alpha_1(\overrightarrow{O_1A_1}^2-r^2)$.
Therefore $c_1(\alpha^1)=0$ if and only if $A_1$ is on $\mathbb{S}^{n-1}_1$.
\end{proof}

To show the geometric properties of $G_{n,2}$ and $G'_{n,2}$,
we give some examples for $n=2$. 
Without loss of generality, assume $\alpha_1>0$ in the following examples.

\begin{figure}[h]
\begin{center}
\resizebox{.7\textwidth}{!}
  {\includegraphics{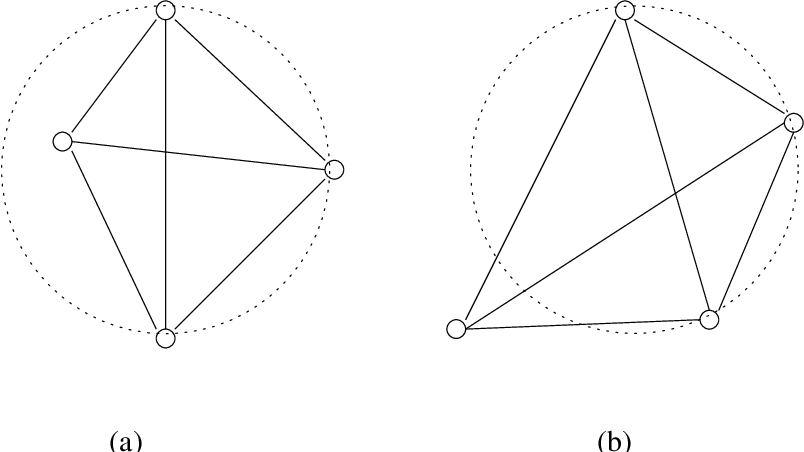}}
\caption{Two 2-tensegrity frameworks in a Euclidean plane in $\mathbb{R}^3$,
with $A_1$ inside or outside of the circle
where $A_2$, $A_3$ and $A_4$ lie on,
which implies $c_1(\alpha^1)<0$ or $c_1(\alpha^1)>0$ respectively}
\label{figure_r_2_2}
\end{center}
\end{figure}

\begin{example}
\label{example_four_points_1}
In Fig. \ref{figure_r_2_2},
assume $A_{1}$, $A_{2}$, $A_{3}$ and $A_{4}$ are the vertices of
a convex quadrilateral in a Euclidean plane in $\mathbb{R}^{3}$.
Topologically, it is hard to tell apart $G_{2,2}$ from $G'_{2,2}$,
because for both frameworks, the quadrilateral is double covered by
volume constraints with opposite signs.
So how to determine that which one of them ($G_{2,2}$ or $G'_{2,2}$)
preserves the volumes of $2$-faces in $\mathbb{R}^{3}$
for small $t\geq 0$?
In Fig. \ref{figure_r_2_2} (a),
$A_1$ is inside the dotted circle
where points $A_2$, $A_3$ and $A_4$ lie on;
and in Fig. \ref{figure_r_2_2} (b) $A_1$ is outside.
Let $O_1$ be the center of the circle and $r$ be the radius, then
$c_1(\alpha^1)=\sum\alpha_i\overrightarrow{O_1A_{i}}^2
=\alpha_1(\overrightarrow{O_1A_1}^2-r^2)$.
As $\alpha_1>0$,
so in (a), $c_1(\alpha^1)<0$ and therefore
by Theorem \ref{theorem_mainone_r_s_h_n_plus_1}
$G'_{2,2}$ preserves the volumes of 2-faces in $\mathbb{R}^{3}$ for small $t\geq 0$;
in (b), $c_1(\alpha^1)>0$ and therefore
$G_{2,2}$ preserves the volumes of 2-faces in $\mathbb{R}^{3}$ for small $t\geq 0$.
\end{example}

\begin{figure}[h]
\begin{center}
\resizebox{.4\textwidth}{!}
  {\includegraphics{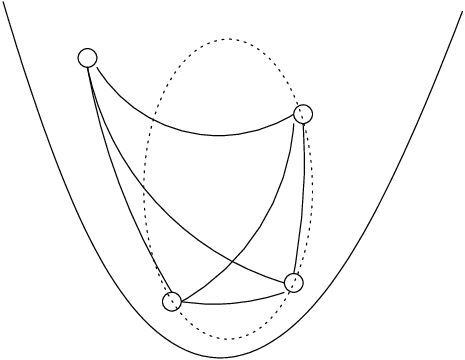}}
\caption{A hyperbolic 2-tensegrity framework in $\mathbb{H}^2$,
where the dotted circle is the intersection 
between $\mathbb{H}^2$ and a $2$-dimensional plane
in $\mathbb{R}^{2,1}$
where $A_2$, $A_3$ and $A_4$ lie on,
with $A_1$ outside of the dotted circle,
which implies $c_1(\alpha^1)>0$}
\label{figure_h_2_2}
\end{center}
\end{figure}

\begin{example}
Fig. \ref{figure_h_2_2} is the hyperbolic version, where the dotted 
circle is the intersection between $\mathbb{H}^2$ and a $2$-dimensional plane in $\mathbb{R}^{2,1}$
where $A_2$, $A_3$ and $A_4$ lie on, with $A_1$ ``outside''
of the dotted circle (but in general the intersection need not be a closed ``circle'').
This implies $A_1$ and the origin $O$ are on the opposite sides of the $2$-dimensional plane.
As $\kappa<0$, $\alpha_1>0$ and 
$c_1(\alpha^1)=\kappa\cdot2\sum\alpha_{i}$ by (\ref{equation_invariant_s_h_A}),
so $c_1(\alpha^1)>0$, and therefore
$G_{2,2}$ preserves the volumes of 2-faces in $\mathbb{H}^3$ for small $t\geq 0$.
\end{example}


From Theorem \ref{theorem_mainone_r_s_h_n_plus_1} and 
Proposition \ref{proposition_common_sphere_hyperplane} we come up with
the following example of ``four points on a circle'', which is rather interesting.

\begin{example}
\label{example_four_points_2}
Still in Fig. \ref{figure_r_2_2},
with four points that are initially in convex position in a 2-dimensional plane in $\mathbb{R}^3$,
but now we constrain all four $2$-faces 
to preserve the volumes while the vertices vary smoothly in $\mathbb{R}^3$.
In order for the four points to form
a \emph{non}-degenerate $3$-simplex in $\mathbb{R}^3$,
by Theorem \ref{theorem_mainone_r_s_h_n_plus_1}
this can only happen when $c_1(\alpha^1)=0$.
Namely, they have to be confined in a plane first
until they move on to a common circle,
and only from this circle they can be \emph{lifted}
to form a non-degenerate $3$-simplex in $\mathbb{R}^3$.
For the same set up in $\mathbb{H}^3$ (Fig. \ref{figure_h_2_2}),
the ``critical position'' for the four points 
to be lifted from $\mathbb{H}^2$ to form a non-degenerate 3-simplex
in $\mathbb{H}^3$ is also when $c_1(\alpha^1)=0$,
namely, when the points are affinely dependent.
The spherical case is also similar, where the ``critical position''
of $c_1(\alpha^1)=0$ is when the four points are on a small circle
in $\mathbb{S}^2$.
\end{example}

\subsection{A positive definite kernel on the hyperbolic space}
\label{section_positive_kernel}

In this section we discuss the positive definiteness of  $g_{\sss F}$
(Definition \ref{definition_partial_derivative_g_theta})
on the hyperbolic space $\mathbb{H}^d$.
Recall that for a set $X$, a symmetric function
$L: X \times X \rightarrow \mathbb{R}$ is a
\emph{positive definite kernel} on $X$
if for any $m\in \mathbb{N}$ and $x_1,\dots,x_m\in X$,
the matrix $(L(x_i,x_j))_{1\leq i,j\leq m}$ is positive semi-definite.
This is also equivalent to having all the principal minors of the matrix
$(L(x_i,x_j))_{1\leq i,j\leq m}$ non-negative.

In $\mathbb{R}^d$, by (\ref{equation_g_B_Euclidean})
$g_{\sss B}(P,Q)=\overrightarrow{PB}\cdot\overrightarrow{QB}$, 
thus $g_{\sss B}$ is a positive definite kernel on $\mathbb{R}^d$.
In $\mathbb{S}^d$ or $\mathbb{H}^d$ of constant curvature $\kappa$,
by Corollary \ref{corollary_g_B},
\begin{equation}
\label{equation_g_B_non_Euclidean}
g_{\sss B}(P,Q)=\frac{2}{1+\kappa P\cdot Q}
\overrightarrow{PB}\cdot\overrightarrow{QB},
\end{equation}
and we have the following analogue for hyperbolic space $\mathbb{H}^1$.

\begin{theorem}
\label{theorem_innerproduct_h}
Let $B\in\mathbb{H}^1$, then
$g_{\sss B}$ is a positive definite kernel on $\mathbb{H}^1$.
\end{theorem}

\begin{proof}
To prove $g_{\sss B}$ is a positive definite kernel on $\mathbb{H}^1$,
it suffices to show that for any $m\in \mathbb{N}$ and $P_1,\dots,P_m\in \mathbb{H}^1$,
\begin{equation}
\label{equation_determinant_positive}
\det(g_{\sss B}(P_i,P_j))_{1\leq i,j\leq m} \geq 0.
\end{equation}

Pick up a direction in $\mathbb{H}^1$, denote the geodesic distance
between $P_i$ and $B$ by $r_i$ if $P_i$ is at the ``right'' side of $B$,
and by $-r_i$ if $P_i$ is at the ``left'' side of $B$.
Then by (\ref{equation_g_B_non_Euclidean}),
\allowdisplaybreaks{
\begin{align*}
&\det(g_{\sss B}(P_i,P_j))_{1\leq i,j\leq m}
=\det\left(\frac{2}{1-P_i\cdot P_j}
\overrightarrow{P_i B}\cdot\overrightarrow{P_j B}\right)  \\
&\quad=\det\left(\frac{2}{1+\cosh(r_i-r_j)}
(\cosh r_i+\cosh r_j- \cosh(r_i-r_j)-1) \right)   \\
&\quad=\det\left(\frac{2}{2\cosh^2\frac{r_i-r_j}{2}}
(2\cosh\frac{r_i+r_j}{2}\cosh\frac{r_i-r_j}{2}
-2\cosh^2\frac{r_i-r_j}{2})\right)                  \\
&\quad=\det\left(\frac{4}{\cosh\frac{r_i-r_j}{2}}
\sinh\frac{r_i}{2}\sinh\frac{r_j}{2} \right)
=2^m\left(\prod_i (e^{r_i}-1) \right)^2
\det\left(\frac{1}{e^{r_i}+e^{r_j}}\right)_{1\leq i,j\leq m}.
\end{align*}
}
By Cauchy's determinant identity, which states that
\[ \det\left(\frac{1}{x_i +y_j}\right)_{1\leq i,j\leq n}
=\frac{\prod_{i<j}(x_j-x_i)(y_j-y_i)}{\prod_{i,j}(x_i+y_j)},
\]
we have
\[\det(g_{\sss B}(P_i,P_j))_{1\leq i,j\leq m}
=2^m\left(\prod_i (e^{r_i}-1) \right)^2
\cdot\frac{\left(\prod_{i<j}(e^{r_j}-e^{r_i})\right)^2}{\prod_{i,j}(e^{r_i}+e^{r_j})}
\geq 0,
\]
which proves (\ref{equation_determinant_positive}) and finishes the proof.
\end{proof}

\begin{remark}
Following a similar proof, we can show that $g_{\sss B}$ is \emph{not} a positive definite kernel
on $\mathbb{S}^1$ for any point $B$ in $\mathbb{S}^1$.
\end{remark}

With Theorem \ref{theorem_innerproduct_h} proved, it seems natural for us to conjecture the following.

\begin{conjecture}
\label{conjecture_innerproduct_h}
Let $B\in\mathbb{H}^d$ and $d\geq 2$,
then $g_{\sss B}$ is a positive definite kernel on $\mathbb{H}^d$.
\end{conjecture}

Next we show that
if $F$ is a $k$-polytope in $\mathbb{R}^d$, 
then $g_{\sss F}$ is a positive definite kernel on $\mathbb{R}^d$.
Let $w_P$ (resp. $w_Q$) be the altitude
vector for the point $P$ (resp. $Q$) with respect to the affine span of $F$,
then it is not hard to show that
$g_{\sss F}(P,Q)=k!V_k(F)w_P\cdot w_Q$,
and the positive definiteness of $g_{\sss F}$ immediately follows.

For the hyperbolic case, we have the following conjecture.

\begin{conjecture}
\label{conjecture_innerproduct_h_F}
Let $F$ be a $k$-polytope in $\mathbb{H}^d$,
then $g_{\sss F}$ is a positive definite kernel on $\mathbb{H}^d$.
\end{conjecture}

\subsection{Related questions and a counterexample}
\label{section_related_questions}

Once Theorem \ref{theorem_mainone_r_s_h_n} and
\ref{theorem_mainone_r_s_h_n_plus_1} are proved,
one question naturally arises:
Under the same condition, is the motion also \emph{rigid}?

This relates to a question raised by Connelly and others:

\begin{question}
\label{question_simplex_congruence}
For $r\geq 2$, do the volumes of all $r$-faces of a
$n$-simplex in $M^n$ determine the $n$-simplex up to congruence?
\end{question}

Question \ref{question_simplex_congruence} was initially posed in the Euclidean space only,
but in the context of this paper, we are also
interested in the spherical and hyperbolic case,
particularly when continuous motion is involved.
The case $r=n-2$ must be classical, 
and various counterexamples were constructed for the Euclidean case
(see \cite{Mc2,MR}).
And following an idea in Mohar and Rivin \cite{MR},
we give a construction for all $r\geq 2$ at the end of this section,
for both Euclidean and non-Euclidean cases.

However, to our knowledge the following continuous analogue of Question
\ref{question_simplex_congruence}
for case $r=n-2$ is still open,
and may very likely to have an affirmative answer.

\begin{question}
\label{question_rigid_motion}
For a $n$-simplex in $M^n$ with $n\geq 4$, 
if a continuous motion preserves the volumes of all
$(n-2)$-faces of the $n$-simplex, then is the motion rigid?
\end{question}

Note that a $n$-simplex
has the same number of edges and $(n-2)$-faces,
and up to congruence is determined by its edge lengths,
so the question is natural.
As the volumes of $(n-2)$-faces are algebraically independent
over the edge lengths (see, for example, \cite{MR}),
Question \ref{question_rigid_motion} should hold
an affirmative answer for ``almost all'' configurations.
While in this paper we do not try to solve Question \ref{question_rigid_motion},
which is mutually independent of our rigidity theorem, an affirmative answer
to Question \ref{question_rigid_motion} will further
improve our main theorems.

\begin{remark}
\label{remark_theorem_mainone}
If Question \ref{question_rigid_motion}
holds an affirmative answer for all non-degenerate simplices,
then Theorem \ref{theorem_mainone_r_s_h_n} and 
\ref{theorem_mainone_r_s_h_n_plus_1}
can be improved to claim that $G_{n,k}$ and $G'_{n,k}$ are rigid for $k\leq n-2$.
If Question \ref{question_rigid_motion}
also holds an affirmative answer for degenerate simplices
with non-degenerate codimension $1$ faces,
then $G_{n,n-1}$ and $G'_{n,n-1}$ are rigid as well.
However, when $n\geq 2$, 
$G_{n,n}$ and $G'_{n,n}$ are never rigid,
as the number of volume constraints is less than the
degree of freedom of $\mathbf{A}$ up to congruence.
\end{remark}

Now we give our construction of a counterexample to 
Question \ref{question_simplex_congruence} for general $r\geq 2$,
essentially following an idea in \cite{MR}.
Let $\Delta_{\epsilon}(t)$ be a $n$-simplex
in $M^n$ whose all sides are equal to a small $\epsilon$ except for one side
whose length is $\epsilon\cdot t$.
It can be shown that for the Euclidean case:
first, $t$ may take any positive value smaller than
$\sqrt{\frac{2n}{n-1}}$; and second,
for any of the $r$-faces that contains the edge with length $\epsilon\cdot t$,
the square of its volume is a quadratic function of $t^2$,
and it
\emph{peaks} when $t=t_0:=\sqrt{\frac{r}{r-1}}$.
As $t_0\leq\sqrt{2}<\sqrt{\frac{2n}{n-1}}$ for $r\geq 2$,
so for sufficiently small $\epsilon$,
$\Delta_{\epsilon}(t_0)$ is \emph{obtainable}
for both Euclidean and non-Euclidean cases.
By properly choosing two close values $t_1$ and $t_2$
satisfying $t_1<t_0<t_2$, 
$\Delta_{\epsilon}(t_1)$ and $\Delta_{\epsilon}(t_2)$ can have the same
volumes on all the corresponding $r$-faces.

\section{Characteristic polynomial of $(\mathbf{A},\alpha)$}
\label{section_char_poly}

For the degenerate $(n+1)$ simplex $\mathbf{A}$
and a $1$-stress $\alpha$ (see (\ref{equation_alpha_r_s_h})),
recall that an invariant $c_{k-1}(\alpha^{k-1})$ (Definition \ref{defintion_invariant_c_k})
plays an important role in a rigidity property of the $k$-faces of $\mathbf{A}$
in Theorem \ref{theorem_mainone_r_s_h_n_plus_1}.
To analyze the relationship between
these rigidity properties of different dimensions $k$,
we introduce a notion of \emph{characteristic polynomial} of
$(\mathbf{A},\alpha)$
by defining
\begin{equation*}
f(x)=\sum_{i=0}^{n+1}(-1)^i c_i(\alpha^i) x^{n+1-i}.
\end{equation*}

Our main result of $f(x)$ is Theorem \ref{theorem_realroots_r},
which shows that the roots of $f(x)$ are real for the Euclidean case,
and gives a way to count the number of positive roots.

\subsection{Properties of the characteristic polynomial}

In this section, let $(\mathbf{A},\alpha)$ be as in (\ref{equation_alpha_r_s_h}),
but for the spherical case
we assume $\mathbf{A}$ is
confined in an open half sphere.
By Remark \ref{remark_invariant_n_plus_1}, we have $c_{n+1}(\alpha^{n+1})=0$.
Let $\{\lambda_{1}, \dots, \lambda_{n}\}$, with no particular order,
be the rest roots of $f(x)$ besides a $0$.

For a $k$-simplex $F$ (as opposed to a more general $k$-polytope)
and two points $P$ and $Q$ in $M^d$, instead of using 
$g_{\sss F}(P,Q)$ (Definition \ref{definition_partial_derivative_g_theta})
sometimes it is more convenient
to use $d_{\sss F}(P,Q)$,
also a new notation introduced in this paper.

\begin{definition}
\label{definition_partial_derivative_d}
For a $k$-simplex $F$ in $M^d$,
define $d_{\sss F}(P,Q)$ by $k!V_k(F)\,g_{\sss F}(P,Q)$ for the Euclidean case,
and by $\norm{F}\,g_{\sss F}(P,Q)$ for the non-Euclidean case.
Also set $d_{\sss\varnothing}(P,Q)=1$.
\end{definition}

\begin{remark}
\label{remark_d_F}
Unlike $g_{\sss F}$ that $F$ need to be non-degenerate, 
$d_{\sss F}$ is well defined when $F$ is degenerate.
\end{remark}

By Definition \ref{defintion_invariant_c_k},
$c_k(\alpha^k)$ can also be equivalently defined by
\begin{equation}
\label{equaition_invariant_r_s_h_A_d}
c_k(\alpha^k) =\sum_{F\subset \mathbf{A},\dim(F)=k-1}
(\prod_{A_s\in F}\alpha_s)d_{\sss F}(P,Q),
\end{equation}
which is independent of the choice of $P$ and $Q$.

In $\mathbb{R}^d$, we give a useful formula for $d_{\sss F}(P,Q)$:
Let $B_1$, \dots, $B_{k+1}$ be the vertices of $F$ in $\mathbb{R}^d$,
then
\begin{equation}
\label{equation_prod_r}
d_{\sss F}(P,Q) =
(\overrightarrow{PB_1}\wedge\cdots\wedge\overrightarrow{PB_{k+1}})
\cdot(\overrightarrow{QB_1}\wedge\cdots\wedge\overrightarrow{QB_{k+1}}).
\end{equation}
When $F$ is degenerate,
the right side of the formula is always well defined,
so we can also extend the definition of $d_{\sss F}$ accordingly.

\begin{proof} [Proof of (\ref{equation_prod_r})]
Denote by $\hat{F}$ the $(k+2)$-dimensional simplex,
which is the join of $F$ with a line segment $PQ$,
and let $\theta_F$ be the dihedral angle at face $F$.
Also let $w_P$ (resp. $w_Q$) be the altitude vector
for point $P$ (resp. $Q$) with respect to the affine span of $F$.
Combine Definition \ref{definition_partial_derivative_d}
and \ref{definition_partial_derivative_g_theta}, we have
\[d_{\sss F}(P,Q) = k!V_k(F)\,
(k+2)!\,\frac{d V_{k+2}(\hat{F})}{d\theta_F}.
\]
Note that
$(k+2)!\,V_{k+2}(\hat{F}) =
k!V_k(F)\,\norm{w_P}\cdot\norm{w_Q}\cdot\sin\theta_F$,
therefore
\begin{align*}
d_{\sss F}(P,Q) 
&=(k!V_k(F))^2\,\norm{w_P}\cdot\norm{w_Q}\cdot\cos\theta_F
=(k!V_k(F))^2\,w_P\cdot w_Q    \\
&=(k!V_k(F)\,w_P)\cdot(k!V_k(F)\,w_Q)             \\
&=(\overrightarrow{PB_1}\wedge\cdots\wedge\overrightarrow{PB_{k+1}})
\cdot(\overrightarrow{QB_1}\wedge\cdots\wedge\overrightarrow{QB_{k+1}}).
\end{align*}
\end{proof}

For the Euclidean case,
without loss of generality, we use the coordinate of $\mathbb{R}^n$
for $\mathbf{A}$ in the following.
Let $B$ be a $(n+1)\times n$ matrix whose $i$-th row is the row vector
$\overrightarrow{A_{n+2}A_i}$ for $i\le n+1$,
and $D=\diag(\alpha_1,\dots,\alpha_{n+1})$
be a diagonal matrix.

\begin{lemma}
The characteristic polynomials of matrix $BB^{T}D$ and $B^{T}DB$
are $f(x)$ and $f(x)/x$ respectively.
\label{lemma_char_poly}
\end{lemma}

\begin{proof}
The coefficient of $x^{n+1-k}$ in the characteristic
polynomial of $BB^{T}D$ is $(-1)^k$ times the sum of
all principal minors of $BB^{T}D$ of order $k$,
which can be shown to be $(-1)^{k}c_k(\alpha^k)$
by choosing $P=Q=A_{n+2}$
in (\ref{equaition_invariant_r_s_h_A_d})
and then applying (\ref{equation_prod_r}).
Therefore $f(x)$ is the
characteristic polynomial of $BB^{T}D$.
Let $B_1$ be $B$ and $B_2$ be $B^{T}D$.
A well known property in linear algebra states that:
If $B_1$ is a $m\times n$ matrix and $B_2$ is a $n\times m$ matrix,
then the characteristic polynomial of $B_1B_2$
is $x^{m-n}$ times the characteristic polynomial of $B_2B_1$.
Therefore $f(x)/x$ is
the characteristic polynomial of $B^{T}DB$.
\end{proof}

Now we have the following main property for $f(x)$.

\begin{theorem}
In the Euclidean case, the roots of $f(x)$ are real.
And if $\{\alpha_{1},\dots,\alpha_{n+2}\}$ has $s$ positive
numbers, then $f(x)$ has $1$ zero,
$s-1$ positive and $n+1-s$ negative roots.
\label{theorem_realroots_r}
\end{theorem}

\begin{proof}
As $\{\lambda_{1},\dots,\lambda_{n}\}$ are the roots of $f(x)/x$,
which by Lemma \ref{lemma_char_poly} is the characteristic polynomial
of a symmetric matrix $B^{T}DB$, so all $\lambda_{i}$ are real. 
In (\ref{equaition_invariant_r_s_h_A_d}) 
by choosing $P=A_1$ and $Q=A_2$ for $k=n$,
we have 
$c_n(\alpha^n)=(\prod_{A_s\in F_{12}}\alpha_s)d_{\sss F_{12}}(A_1,A_2)$,
where $F_{12}$ is the $(n-1)$-face of $\mathbf{A}$ 
that without the vertices $A_1$ and $A_2$.
Then by (\ref{equation_prod_r}) we have $c_n(\alpha^n)\neq 0$, 
and thus all $\lambda_{i}$ are also non-zero.
So $\{\lambda_{1},\dots,\lambda_{n}\}$ must have the same signs as a $n$-subset
of the diagonal entries of $D=\diag(\alpha_1,\dots,\alpha_{n+1})$.
By symmetry, $\{\lambda_{1},\dots,\lambda_{n}\}$ should have
the same signs as a $n$-subset of
$\{\alpha_1,\dots,\hat{\alpha_j},\dots,\alpha_{n+2}\}$
for any $j$ with $1\leq j\leq n+2$.
So $\{\lambda_{1},\dots,\lambda_{n}\}$ must have $s-1$ positive
and $n+1-s$ negative numbers.
\end{proof}

\begin{theorem}
In the Euclidean case,
$f(x)/x$ has $n$-repeated roots if and only if
$\overrightarrow{A_{i}A_{j}}\cdot\overrightarrow{A_{k}A_{l}}=0$
for all distinct numbers $i$, $j$, $k$ and $l$.
\label{theorem_repeatedroots}
\end{theorem}

\begin{proof}
The ``only if'' part. If $\lambda_{1}=\cdots=\lambda_{n}$, denote it by $\lambda$.
Then by Lemma \ref{lemma_char_poly}, $B^{T}DB=\lambda I_n$ where $I_n$
is the $n\times n$ identity matrix.
Let $B_1$ be a $(n+1)\times (n+1)$ matrix, whose first $n$ columns are $B$,
and every entry on the last column is $\sqrt{-\lambda/\alpha_{n+2}}$
(it is ok if it is not real).
Easy to see that $B_1^{T}DB_1=\lambda I_{n+1}$.
So $B_1^T=\lambda B_1^{-1}D^{-1}$
and therefore $B_{1}B_1^{T}=\lambda D^{-1}$.
So if $i\neq j$,
then the $(i,j)$-th entry of $B_{1}B_1^{T}$ is $0$,
and therefore
\begin{equation}
\label{equation_lambda_equal}
\overrightarrow{A_{n+2}A_i}\cdot\overrightarrow{A_{n+2}A_j}
-\frac{\lambda}{\alpha_{n+2}}=0.
\end{equation}
If $i$, $j$, $k$ and $n+2$ are distinct, replace ``$j$''
with ``$k$'' in (\ref{equation_lambda_equal}) and subtract
it from (\ref{equation_lambda_equal}),
then
$\overrightarrow{A_{n+2}A_i}\cdot\overrightarrow{A_jA_k}=0$.
By symmetry,
$\overrightarrow{A_{i}A_{j}}\cdot\overrightarrow{A_{k}A_{l}}=0$
for all distinct $i$, $j$, $k$ and $l$.

The ``if'' part. Assume
$\overrightarrow{A_{i}A_{j}}\cdot\overrightarrow{A_{k}A_{l}}=0$
for all distinct $i$, $j$, $k$ and $l$.
Then 
\begin{equation}
\overrightarrow{A_{i}A_{j}}\cdot\overrightarrow{A_{i}A_{k}}
=\overrightarrow{A_{i}A_{j}}\cdot\overrightarrow{A_{i}A_{l}}.
\end{equation}
So for a fixed $i$,
$\overrightarrow{A_{i}A_{j}}\cdot\overrightarrow{A_{i}A_{k}}$
is independent of $j$ and $k$, as long as $i$, $j$ and $k$ are distinct.
We denote it by $b_i$.
Let $B_2$ be a $n\times n$ matrix whose $i$-th row is vector
$\overrightarrow{A_{n+1}A_i}$ for $i\le n$, 
$B_3$ be a $n\times n$
matrix whose $i$-th row is vector
$\overrightarrow{A_{n+2}A_i}$,
and $D_1=\diag(\alpha_1,\dots,\alpha_n)$ be a diagonal matrix.
By choosing $P=A_{n+1}$ and $Q=A_{n+2}$ in
(\ref{equaition_invariant_r_s_h_A_d})
and then applying (\ref{equation_prod_r}),
$f(x)/x$ is the characteristic polynomial
of $B_{2}B_{3}^{T}D_1$. Since
$\overrightarrow{A_{n+1}A_i}\cdot\overrightarrow{A_{n+2}A_j}$
is $b_i$ when $i=j$ and 0 when $i\ne j$, then
$B_{2}B_{3}^{T}D_1=\diag(\alpha_i b_i)_{i\le n}$.
So $\{\lambda_1,\dots,\lambda_n\}$ is $\{\alpha_i b_i\}_{i\le n}$
in some order.
By symmetry, any $n$-subset of $\{\alpha_i b_i\}_{i\leq n+2}$
is $\{\lambda_1,\dots,\lambda_n\}$ in some order as well.
So all $\alpha_i b_i$ are equal, and therefore $\lambda_1=\dots=\lambda_n$.
\end{proof}

It is natural to ask if $f(x)$ still has real roots in the
non-Euclidean case,
and some evidence suggests
that the hyperbolic version of
Theorem \ref{theorem_realroots_r} might still hold.
We numerically computed some examples for $n=2$, 
and our tests for the hyperbolic space all have real roots,
but the same test shows that the spherical version of
Theorem \ref{theorem_realroots_r} is not true.

\begin{conjecture}
\label{conjecture_realroots_h}
In the hyperbolic case, the roots of $f(x)$ are real.
And if $\{\alpha_{1},\dots,\alpha_{n+2}\}$ has $s$ positive numbers,
then $f(x)$ has $1$ zero, $s-1$ positive and $n+1-s$ negative roots.
\end{conjecture}

\subsection{Generalization of characteristic polynomial}
\label{section_char_poly_general}

So far the characteristic polynomial $f(x)$ is defined on degenerate
$(n+1)$-simplices only.
To complete the discussion of $f(x)$,
with the proofs skipped,
we loosely discuss a generalization of $f(x)$
by showing that it can be naturally
generalized to a general $(\mathbf{A}, \alpha)$ (not necessarily finite) in $M^n$. 
We start with a finite set.

Abuse of notation:
Let $\mathbf{A}=\{A_1, \dots, A_m\}$ be a set of $m$ $(m\geq n+2)$ points in $M^n$
in general position,
and $\alpha=\{\alpha_{1}, \dots, \alpha_m\}$ be a $1$-stress on $\mathbf{A}$.
For this \emph{new} $(\mathbf{A},\alpha)$ and each $k\leq n+1$,
by Definition \ref{defintion_invariant_c_k} and Theorem \ref{theorem_invariant_r_s_h}
we have
\begin{equation}
\label{equaition_invariant_r_s_h_A_d_new}
\sum_{F\subset \mathbf{A},\dim(F)=k-1}
(\prod_{A_s\in F}\alpha_s)d_{\sss F}(P,Q) = c_k(\alpha^k),
\end{equation}
where the $k$-stress $\alpha^k$ is by Definition \ref{definition_k_stress_alpha}
and $c_k(\alpha^k)$ is an invariant
independent of the choice of $P$ and $Q$;
and for the non-Euclidean case,
\begin{equation}
\label{equation_invariant_s_h_A_copy}
c_k(\alpha^k) =\kappa (k+1)(k-1)!
\sum_{F\subset \mathbf{A},\dim(F)=k-1}(\prod_{A_s\in F}\alpha_s)
\norm{F}
\,V_{k-1}(F).
\end{equation}

\begin{remark}
The right side of (\ref{equation_invariant_s_h_A_copy})
is well defined when $F$ is allowed to be degenerate. 
Even if the term $V_{k-1}(F)$ is
not well defined for some reason, say, $F$ contains a pair of antipodal
points in the spherical case, then the term $\norm{F}$ is always zero
and will make their product zero.
This suggests the possibility to define $c_k(\alpha^k)$ in a much more general sense,
say, when $(\mathbf{A}, \alpha)$ is distributed in a continuous manner.
\end{remark}


Define $f(x)=\sum_{i=0}^{n+1}(-1)^i c_i(\alpha^i) x^{n+1-i}$
as the characteristic polynomial of $(\mathbf{A}, \alpha)$ as before.
We want to point out that $c_{n+1}(\alpha^{n+1})$ still vanishes unless
$\mathbf{A}$ is not confined in any open half sphere in the spherical case.

Slightly modifying Lemma \ref{lemma_char_poly}
and Theorem \ref{theorem_realroots_r}, we can prove
that the roots of this newly generalized
$f(x)$ are still real in the Euclidean case.
And again, we conjecture
the same for the hyperbolic case.

\begin{conjecture}
In the hyperbolic case, the roots of the generalized $f(x)$ are real.
\end{conjecture}

\noindent
{\bf Acknowledgements:}
This article is an extension of the author's Ph.D. thesis \cite{Zh}
at M.I.T.. I am very grateful to my advisor Professor D. Kleitman
and to Professor R. Stanley for their guidance during this work.
I would also like to thank Professor
R. Connelly, Wei Luo and Xun Dong
for their many helpful suggestions and discussions.


{\footnotesize


}

\end{document}